\documentclass[11pt]{article}
\usepackage{amsmath}
\usepackage{amssymb}

\usepackage{theorem}
\numberwithin{equation}{section}

\newtheorem{theorem}{Theorem}[section]

\newtheorem{corollary}[theorem]{Corollary}
\newtheorem{lemma}[theorem]{Lemma}
\newtheorem{definition}{Definition}[section]

\begin{document}
\title{Global well - posedness and scattering for nonlinear Schr{\"o}dinger equations with algebraic nonlinearity when $d = 2, 3$ and $u_{0}$ is radial}
\date{\today}
\author{Benjamin Dodson}
\maketitle

\noindent \textbf{Abstract:} In this paper we discuss global well-posedness and scattering for some initial value problems that are $\dot{H}^{1}$ subcritical. We prove global well-posedness and scattering for radial data in $H^{s}$, $s > s_{c}$, where the initial value problem is $\dot{H}^{s_{c}}$-critical. We make use of the long time Strichartz estimates of \cite{D2} to do this.

\section{Introduction}
In this paper we examine the three dimensional initial value problem

\begin{equation}\label{1.1}
(i \partial_{t} + \Delta) u = F(u) = |u|^{2} u, \hspace{5mm} u(0,x) = u_{0} \in H_{x}^{s}(\mathbf{R}^{3}),
\end{equation}

\noindent as well as the two dimensional initial value problems

\begin{equation}\label{1.2}
(i \partial_{t} + \Delta) u = |u|^{k} u, \hspace{5mm} u(0,x) = u_{0} \in H_{x}^{s}(\mathbf{R}^{2}),
\end{equation}

\noindent where $k$ may be any positive integer. In each case $u_{0}$ is a radial function.\vspace{5mm}

\noindent Solutions to $(\ref{1.1})$ and $(\ref{1.2})$ give rise to a family of solutions via the scaling,

\begin{equation}\label{1.3}
 u(t,x) \mapsto u_{\lambda}(t,x) = \lambda^{\frac{1}{k}} u(\lambda^{2} t, \lambda x).
\end{equation}

\noindent Under this scaling, for any $s \in \mathbf{R}$,

\begin{equation}\label{1.4}
 \| u_{\lambda}(0,x) \|_{\dot{H}_{x}^{s}(\mathbf{R}^{d})} = \lambda^{-\frac{d}{2} + s + \frac{1}{k}} \| u(0, x) \|_{\dot{H}_{x}^{s}(\mathbf{R}^{d})}.
\end{equation}

\noindent Thus, $(\ref{1.1})$ is called $\dot{H}^{1/2}$-critical since under $(\ref{1.4})$,

\begin{equation}\label{1.5}
\| u_{\lambda}(0, x) \|_{\dot{H}_{x}^{1/2}(\mathbf{R}^{3})} = \| u(0, x) \|_{\dot{H}_{x}^{1/2}(\mathbf{R}^{3})}.
\end{equation}

\noindent Likewise, $(\ref{1.2})$ is called $\dot{H}^{s_{c}}$-critical, where $s_{c} = 1 - \frac{1}{k}$. 

This scaling is crucial to local well-posedness. Recall the usual definition of well-posedness.

\begin{definition}[Well-posedness]\label{d1.2}
The initial value problem $(\ref{1.1})$ is well-posed on an open interval $I \subset \mathbf{R}$, $0 \in I$, for $u_{0} \in H_{x}^{s}(\mathbf{R}^{3})$, if

\begin{enumerate}
\item $(\ref{1.1})$ has a unique solution $u$ lying in $C_{t}^{0}(I ; H_{x}^{s}(\mathbf{R}^{3}))$,

\item The solution satisfies the Duhamel formula

\begin{equation}\label{1.8}
u(t) = e^{it \Delta} u_{0} - i \int_{0}^{t} e^{i(t - \tau) \Delta} (|u|^{2} u)(\tau) d\tau,
\end{equation}

\item For any compact $J \subset I$, the map $u_{0} \mapsto L_{t,x}^{5}(J \times \mathbf{R}^{3})$ is continuous.
\end{enumerate}

\noindent The definition of global well - posedness for $(\ref{1.2})$ corresponds to $(1)-(3)$ above, although $L_{t,x}^{5}$ should be replaced by $L_{t,x}^{4k}$ and $\mathbf{R}^{3}$ should be replaced by $\mathbf{R}^{2}$.

\end{definition}

Also recall the definition of scattering.

\begin{definition}[Scattering]\label{d1.3}
A global solution to $(\ref{1.1})$ and $(\ref{1.2})$ with initial data $u_{0}$ is said to scatter forward in time to some $u_{+} \in H_{x}^{s}(\mathbf{R}^{d})$ if

\begin{equation}\label{1.9}
 \lim_{t \rightarrow +\infty} \| u(t) - e^{it \Delta} u_{+} \|_{H^{s}(\mathbf{R}^{d})} = 0.
\end{equation}

\noindent Analogously, $u$ is said to scatter backward in time to some $u_{-} \in H_{x}^{s}(\mathbf{R}^{3})$ if

\begin{equation}\label{1.9.1}
 \lim_{t \rightarrow -\infty} \| u(t) - e^{it \Delta} u_{-} \|_{H^{s}(\mathbf{R}^{d})} = 0.
\end{equation}

\noindent $(\ref{1.1})$ is said to be scattering for initial data lying in a certain set $X$ if for each $u_{0} \in X$ there exists $u_{+}$ and $u_{-}$ such that $(\ref{1.9})$ and $(\ref{1.9.1})$ hold, and furthermore, the maps $u_{0} \mapsto u_{+}$ and $u_{0} \mapsto u_{-}$ are continuous as functions of $u_{0}$.
\end{definition}

\noindent \textbf{Remark:} Scattering for $(\ref{1.1})$ corresponds to $\| u \|_{L_{t,x}^{5}(\mathbf{R} \times \mathbf{R}^{3})} < \infty$ and scattering for $(\ref{1.2})$ corresponds to $\| u \|_{L_{t,x}^{4k}(\mathbf{R} \times \mathbf{R}^{2})} < \infty$.

\begin{theorem}\label{t1.1}
$(\ref{1.1})$ is locally well-posed for any $u_{0} \in H^{s}(\mathbf{R}^{3})$, $s > \frac{1}{2}$ on some interval $[-T, T]$, $T(\| u_{0} \|_{H^{s}}, s) > 0$. If $u_{0} \in \dot{H}^{1/2}(\mathbf{R}^{3})$ then $(1.1)$ is locally well-posed on some interval $[-T, T]$, $T(u_{0}) > 0$, where $T(u_{0})$ depends on the profile of the initial data and not just its size. Moreover, for $\| u_{0} \|_{\dot{H}^{1/2}(\mathbf{R}^{3})}$ small, $(\ref{1.1})$ is globally well-posed and scattering.\vspace{5mm}

\noindent The corresponding results also hold for $(\ref{1.2})$ and the critical space $\dot{H}^{1 - \frac{1}{k}}(\mathbf{R}^{2})$.
\end{theorem}

\noindent \emph{Proof:} See \cite{CaWe}, \cite{CaWe1}. $\Box$\vspace{5mm}

\noindent \textbf{Remark:} \cite{ChrColTao2} and \cite{ChrColTao1} proved that Theorem $\ref{t1.1}$ is sharp.\vspace{5mm}

\noindent \textbf{Remark:} \cite{KM2} proved that $(\ref{1.1})$ is globally well-posed and scattering if and only if $\| u(t) \|_{\dot{H}_{x}^{1/2}(\mathbf{R}^{3})}$ is uniformly bounded on its interval of existence. See \cite{XYu} for the same result when $k = 1$ and $d = 2$.\vspace{5mm}

\noindent $(\ref{1.2})$ with $k = 1$ is now completely solved. \cite{KTV} proved that $(\ref{1.2})$ is globally well-posed and scattering for any $u_{0} \in L^{2}(\mathbf{R}^{2})$, $u_{0}$ radial. \cite{D4} extended this to nonradial data.\vspace{5mm}

\noindent In this paper we show that $(\ref{1.1})$ and $(\ref{1.2})$ are globally well-posed and scattering for $u_{0} \in H_{x}^{s}(\mathbf{R}^{d})$, $u_{0}$ radial, $s > \frac{1}{2}$, and $s > 1 - \frac{1}{k}$ respectively. 

We begin with the cubic problem in three dimensions.
\begin{theorem}\label{t1.4}
The initial value problem $(\ref{1.1})$ is globally well-posed and scattering for any $s > \frac{1}{2}$ for $u_{0}$ radial.
\end{theorem}

\noindent Next, we will prove an explicit upper bound on the scattering size, or $L_{t,x}^{4}$ norm, for a solution to the two dimensional, cubic problem ($k = 1$) for $u_{0}$ radial lying in a subspace of $L_{x}^{2}(\mathbf{R}^{2})$. Scattering for the two dimensional, radial, cubic problem has already been proved for \cite{KTV}. See \cite{D4} for a proof in the nonradial case. However, no explicit norm was computed in \cite{KTV} or \cite{D4}, which we will do here.

\begin{theorem}\label{t1.4.1}
When $k = 1$ and $u_{0}$ is radial, $(\ref{1.2})$ has a global solution with

\begin{equation}\label{1.9.2}
\| u \|_{L_{t,x}^{4}(\mathbf{R} \times \mathbf{R}^{2})}^{4} \lesssim (\| u_{0} \|_{\dot{H}^{s}(\mathbf{R}^{2})} + \| |x|^{s} u_{0} \|_{L^{2}(\mathbf{R}^{2})})^{\frac{8(1 - s)}{s} + 2} (1 + \| u_{0} \|_{L^{2}})^{\frac{4(1 - s)}{s} + 1}.
\end{equation}

\noindent A solution to the focusing problem

\begin{equation}
i u_{t} + \Delta u = -|u|^{2} u, \hspace{5mm} u(0, x) = u_{0},
\end{equation}

\noindent has the scattering size bound

\begin{equation}\label{1.9.3}
\| u \|_{L_{t,x}^{4}(\mathbf{R} \times \mathbf{R}^{2})}^{4} \lesssim (\| u_{0} \|_{\dot{H}^{s}(\mathbf{R}^{2})} + \| |x|^{s} u_{0} \|_{L^{2}(\mathbf{R}^{2})})^{\frac{8(1 - s)}{s} + 2} (1 + \| u_{0} \|_{L^{2}})^{\frac{4(1 - s)}{s} + 1} (1 - \frac{\| u_{0} \|_{L^{2}}^{2}}{\| Q \|_{L^{2}}^{2}})^{-\frac{1}{s}}.
\end{equation}

\noindent $Q$ is the ground state of the focusing problem, that is, the positive solution to

\begin{equation}
\Delta Q + Q^{3} = Q.
\end{equation}
\end{theorem}
\textbf{Remark:} \cite{MarTay} proved a result in this form for $s = 1$.\medskip

\noindent Finally we prove two dimensional scattering results for $(\ref{1.2})$ when $k > 1$.

\begin{theorem}\label{t1.5}
The initial value problem $(\ref{1.2})$ is globally well-posed and scattering for any $s > 1 - \frac{1}{k}$, $u_{0}$ radial. 
\end{theorem}

\subsection{Method of proof} The I-method is used to prove Theorems $\ref{t1.4}$, $\ref{t1.4.1}$, and $\ref{t1.5}$. A solution to $(\ref{1.1})$ conserves the quantities mass,

\begin{equation}\label{1.10}
M(u(t)) = \int |u(t,x)|^{2} dx = M(u(0)),
\end{equation}

\noindent and energy,

\begin{equation}\label{1.11}
E(u(t)) = \frac{1}{2} \int |\nabla u(t,x)|^{2} dx + \frac{1}{4} \int |u(t,x)|^{4} dx.
\end{equation}

\noindent Likewise, a solution to $(\ref{1.2})$ conserves mass $(\ref{1.10})$ and the energy

\begin{equation}\label{1.12}
E(u(t)) = \frac{1}{2} \int |\nabla u(t,x)|^{2} dx + \frac{1}{2k + 2} \int |u(t,x)|^{2k + 2} dx.
\end{equation}

\noindent $(\ref{1.11})$ and $(\ref{1.12})$ combined with the local well-posedness theorem (Theorem $\ref{t1.1}$) proves that $(\ref{1.1})$ and $(\ref{1.2})$ are globally well-posed for data in $H^{1}$. See \cite{GV1}, \cite{LinStrauss} for a proof of scattering in the radial case; \cite{CGT1}, \cite{CKSTT2}, and \cite{PV} for a proof of scattering in the nonradial case for $u_{0} \in H^{1}$.\vspace{5mm}

\noindent The reason for the gap between the local well-posedness result of Theorem $\ref{t1.1}$ and the regularity needed to prove a global result in \cite{GV1} ($s = 1$) is due to an absence of a conserved quantity that controls $\| u(t) \|_{\dot{H}^{s}}$ for $0 < s < 1$. It is true that the momentum, a $\dot{H}^{1/2}$-critical quantity, is conserved, but this quantity does not control the $\dot{H}^{1/2}$ norm.\vspace{5mm}

\noindent The first progress in extending the global well-posedness results for data in $H^{1}$ to $H^{s}$, $s < 1$ came from the Fourier truncation method. \cite{B2} proved that the cubic nonlinear initial value problem is globally well-posed in two dimensions for data in $H^{s}$, $s > \frac{3}{5}$ when $d = 2$. In three dimensions \cite{B6} proved global well-posedness for $s > \frac{11}{13}$ and global well-posedness and scattering for $s > \frac{5}{7}$ when $u_{0}$ is radial. In fact, \cite{B2}, \cite{B6} proved something more, namely that for $s$ in the appropriate interval

\begin{equation}\label{1.13}
u(t) - e^{it \Delta} u_{0} \in H^{1}(\mathbf{R}^{d}).
\end{equation}

\noindent It was precisely $(\ref{1.13})$ that lead to the development of the I-method since $(\ref{1.13})$ is false for many dispersive partial differential equations. See \cite{KTwm} for example. Instead, \cite{CKSTT1} defined an operator $I : H^{s}(\mathbf{R}^{d}) \rightarrow H^{1}(\mathbf{R}^{d})$. Tracking the change of $E(Iu(t))$, \cite{CKSTT1} proved global well-posedness for the cubic nonlinear Schr{\"o}dinger equation when $d = 2$ for $s > \frac{4}{7}$, and when $d = 3$ for $s > \frac{5}{6}$. \cite{CKSTT2} extended the $d = 3$ result to $s > \frac{5}{6}$. $\cite{D1}$ extended this to $s > \frac{5}{7}$, and then \cite{Su} extended this result to $s > \frac{2}{3}$.\vspace{5mm}

\noindent Both \cite{D1} and \cite{Su} utilized the linear-nonlinear decomposition. See also \cite{Roy} for this method in the context of the wave equation. Here we will use the long time Strichartz estimates of \cite{D2}. We show that for radial data, the long time Strichartz estimates decay rapidly, and thus can beat any polynomial power of $N$ arising from the $I$-operator.\vspace{5mm}

\noindent The outline of the paper is as follows. In $\S 2$ we will recall some linear estimates needed in the proof. In $\S 3$ we will describe the I-method and outline the proof of Theorems $\ref{t1.4}$ and $\ref{t1.5}$. In $\S 4$ we will make an induction on frequency argument and prove long time Strichartz estimates for $d = 3$. In $\S 5$ we will prove the energy increment in $d = 3$, yielding Theorem $\ref{t1.4}$. In $\S 6$ we prove Theorem $\ref{t1.4.1}$, obtaining scattering size for the cubic problem in dimension $d = 2$. Then in $\S 7$ we will make an induction on frequency argument and prove long-time Strichartz estimates for $d = 2$ when $k > 1$. Finally in $\S 8$ we will prove the energy increment in $d = 2$, yielding Theorem $\ref{t1.5}$.\vspace{5mm}

\noindent At this point it is necessary to mention some notation used in the paper. This notation was used in \cite{CKSTT1}. The expression $A \lesssim_{B} D$ indicates $A \leq C(B) D$, where $C(B)$ is some constant. When we say $A \lesssim_{\| u_{0} \|_{H^{s}}} B$ or $A \lesssim_{\| u_{0} \|_{H^{s}}, k} B$ we mean that $A \leq C(\| u_{0} \|_{H^{s}}, s, k) B$. $A \sim B$ denotes $A \lesssim B$ and $B \lesssim A$.\vspace{5mm}

\noindent We will also use the notation $A \lesssim B^{a+}$. This means that for any $\epsilon > 0$, there exists $C(\epsilon)$ such that $A \leq C(\epsilon) B^{a + \epsilon}$. We will also use expressions like $\| u \|_{L^{p+}} \lesssim A$, which means that $\| u \|_{L^{p + \epsilon}} \leq C(\epsilon) A$. $\| u \|_{L^{p+}} \gtrsim A$ has the obvious definition.\vspace{5mm}

\noindent Throughout the paper it is unnecessary to distinguish between $u$ and $\bar{u}$. Therefore, we will often write expressions like $|u|^{2} u$ as $u^{3}$ for convenience.\medskip

\noindent \textbf{Acknowledgements:} During the time of writing this paper, the author was supported by an NSF postdoc, as well as NSF grants DMS-$1500424$ and DMS-$1764358$. The author was also a member of the Institute for Advanced Study as a von Neumann fellow during part of the writing of this paper. Finally, the author would like to acknowledge the helpful suggestions of the anonymous referee.

\section{Linear estimates}
In this section we mention a number of estimates for the linear Schr{\"o}dinger equation. None of the results in this section are new.

\subsection{Sobolev spaces}
\begin{definition}[Littlewood-Paley decomposition]\label{d7.1}
Take $\psi \in C_{0}^{\infty}(\mathbf{R}^{d})$, $\psi(x) = 1$ for $|x| \leq 1$, $\psi = 0$ for $|x| > 2$, where $\psi(x)$ is radial and decreasing. Then for any $j$ let

\begin{equation}\label{7.1}
\phi_{j}(x) = \psi(2^{-j} x) - \psi(2^{-j + 1} x).
\end{equation}

\noindent Let $P_{j}$ be the Fourier multiplier given by

\begin{equation}\label{7.2}
\widehat{P_{j} f}(\xi) = \phi_{j}(\xi) \hat{f}(\xi).
\end{equation}

\noindent This gives the Littlewood-Paley decomposition

\begin{equation}
f = \sum_{j = -\infty}^{\infty} P_{j} f,
\end{equation}

\noindent at least in the $L^{2}$ sense.
\end{definition}

\noindent The Littlewood-Paley decomposition is quite useful since

\begin{theorem}[Littlewood-Paley theorem]\label{t7.2}
For any $1 < p < \infty$,

\begin{equation}\label{7.3}
\| f \|_{L^{p}(\mathbf{R}^{d})} \sim_{p, d} \| (\sum_{j = -\infty}^{\infty} |P_{j} f|^{2})^{1/2} \|_{L^{p}(\mathbf{R}^{d})}.
\end{equation}
\end{theorem}

\begin{definition}[Sobolev spaces]\label{d7.3}
For $s \in \mathbf{R}$ the Sobolev space $\dot{H}^{s}(\mathbf{R}^{d})$ is the space of functions whose Fourier transform has finite weighted $L^{2}$ norm,

\begin{equation}\label{7.4}
\| f \|_{\dot{H}^{s}(\mathbf{R}^{d})} = \| |\xi|^{s} \hat{f}(\xi) \|_{L^{2}(\mathbf{R}^{d})},
\end{equation}

\noindent where

\begin{equation}\label{7.5}
\hat{f}(\xi) = (2 \pi)^{-d/2} \int e^{-ix \cdot \xi} f(x) dx.
\end{equation}

\noindent We define the inhomogeneous space

\begin{equation}\label{7.6}
\| f \|_{H^{s}(\mathbf{R}^{d})} = \| (1 + |\xi|^{2})^{s/2} \hat{f}(\xi) \|_{L^{2}(\mathbf{R}^{d})}.
\end{equation}

\noindent Notice that

\begin{equation}\label{7.7}
\| P_{j} f \|_{L^{2}(\mathbf{R}^{d})} \lesssim 2^{-js} \| f \|_{\dot{H}^{s}(\mathbf{R}^{d})}, \hspace{5mm} \| P_{j} f \|_{L^{2}(\mathbf{R}^{d})} \lesssim \inf(2^{-js}, 1) \| f \|_{H^{s}(\mathbf{R}^{d})}.
\end{equation}
\end{definition}

\noindent \textbf{Remark:} $(\ref{7.7})$ is called Bernstein's inequality.\vspace{5mm}

\noindent It follows from H{\"o}lder's inequality that for $2 \leq p \leq \infty$,

\begin{equation}\label{7.8}
\| P_{j} f \|_{L^{p}(\mathbf{R}^{d})} \lesssim_{d} 2^{jd(\frac{1}{2} - \frac{1}{p})} \| P_{j} f \|_{L^{2}(\mathbf{R}^{d})}.
\end{equation}

\noindent Then for $1 < p < \infty$, $s = d(\frac{1}{2} - \frac{1}{p})$,

\begin{equation}\label{7.9}
\| f \|_{L^{p}(\mathbf{R}^{d})} \lesssim_{s, d} \| f \|_{\dot{H}^{s}(\mathbf{R}^{d})}.
\end{equation}

\noindent We also have the radial Sobolev embedding

\begin{equation}\label{7.10}
\| |x| P_{j} f \|_{L^{\infty}(\mathbf{R}^{3})} \lesssim \| P_{j} f \|_{\dot{H}^{1/2}(\mathbf{R}^{3})}.
\end{equation}

\noindent See \cite{St1}, \cite{St}, \cite{T3}, \cite{T1}, and many other sources for more details on Sobolev spaces.

\subsection{Strichartz estimates}
\begin{theorem}\label{t7.4}
Let $e^{it \Delta}$ be the solution operator to the linear evolution equation $(i \partial_{t} + \Delta) u = 0$. That is, $u = e^{it \Delta} u_{0}$ solves

\begin{equation}\label{7.11}
(i \partial_{t} + \Delta)u = 0, \hspace{5mm} u(0,x) = u_{0}.
\end{equation}

\noindent When $d = 3$ define

\begin{equation}\label{7.12}
(p, q) \in \mathcal A_{3} \hspace{5mm} \text{if and only if} \hspace{5mm} 2 \leq p \leq \infty, \hspace{5mm} \text{and} \hspace{5mm} \frac{2}{p} = 3(\frac{1}{2} - \frac{1}{p}).
\end{equation}

\noindent When $d = 2$ define

\begin{equation}\label{7.13}
(p, q) \in \mathcal A_{2} \hspace{5mm} \text{if and only if} \hspace{5mm} 2 < p \leq \infty, \hspace{5mm} \text{and} \hspace{5mm} \frac{1}{p} + \frac{1}{q} = \frac{1}{2}.
\end{equation}

\noindent If $(p, q)$ lies in $\mathcal A_{d}$ then we say that $(p, q)$ is an admissible pair.\vspace{5mm}

\noindent Let $p'$ denote the Lebesgue dual, that is $\frac{1}{p} + \frac{1}{p'} = 1$. Then if $(p, q)$ and $(\tilde{p}, \tilde{q})$ are admissible pairs, when $d = 2$,

\begin{equation}\label{7.14}
\aligned
\| e^{it \Delta} u_{0} \|_{L_{t}^{p} L_{x}^{q}(\mathbf{R} \times \mathbf{R}^{2})} \lesssim_{p} \| u_{0} \|_{L^{2}(\mathbf{R}^{2})}, \\
\| \int_{0}^{t} e^{i(t - \tau) \Delta} F(\tau) d\tau \|_{L_{t}^{p} L_{x}^{q}(I \times \mathbf{R}^{2})} \lesssim_{p, \tilde{p}} \| F \|_{L_{t}^{\tilde{p}'} L_{x}^{\tilde{q}'}(I \times \mathbf{R}^{2})},
\endaligned
\end{equation}

\noindent and when $d = 3$,

\begin{equation}\label{7.15}
\aligned
\| e^{it \Delta} u_{0} \|_{L_{t}^{p} L_{x}^{q}(\mathbf{R} \times \mathbf{R}^{3})} \lesssim \| u_{0} \|_{L^{2}(\mathbf{R}^{3})}, \\
\| \int_{0}^{t} e^{i(t - \tau) \Delta} F(\tau) d\tau \|_{L_{t}^{p} L_{x}^{q}(I \times \mathbf{R}^{3})} \lesssim \| F \|_{L_{t}^{\tilde{p}'} L_{x}^{\tilde{q}'}(I \times \mathbf{R}^{3})},
\endaligned
\end{equation}
\end{theorem}

\noindent \emph{Proof:} \cite{Stri} proved this theorem in the case $p = q$, $\tilde{p} = \tilde{q}$. See \cite{ChrKis}, \cite{GV}, and \cite{Yaj} for a proof of the general result, $p > 2$. \cite{KT} proved the endpoint result $p = 2$ when $d = 3$. \cite{Tao} gives a nice description of the overall theory. $\Box$\vspace{5mm}

\noindent Because of this fact it is convenient, especially in three dimensions, to work with the Strichartz space and the dual Strichartz space.

\begin{definition}[Strichartz space]\label{d7.5}
Let $S^{0}$ be the Strichartz space

\begin{equation}\label{7.16}
S^{0}(I \times \mathbf{R}^{3}) = L_{t}^{\infty} L_{x}^{2}(I \times \mathbf{R}^{3}) \cap L_{t}^{2} L_{x}^{6}(I \times \mathbf{R}^{3}).
\end{equation}

\noindent Let $N^{0}$ be the dual

\begin{equation}\label{7.17}
N^{0}(I \times \mathbf{R}^{3}) = L_{t}^{1} L_{x}^{2}(I \times \mathbf{R}^{3}) + L_{t}^{2} L_{x}^{6/5}(I \times \mathbf{R}^{3}).
\end{equation}

\noindent Then Theorem $\ref{t7.4}$ implies

\begin{equation}\label{7.18}
\aligned
\| e^{it \Delta} u_{0} \|_{S^{0}(\mathbf{R} \times \mathbf{R}^{3})} \lesssim \| u_{0} \|_{L^{2}(\mathbf{R}^{3})}, \\
\| \int_{0}^{t} e^{i(t - \tau) \Delta} F(\tau) d\tau \|_{S^{0}(I \times \mathbf{R}^{3})} \lesssim \| F \|_{N^{0}(I \times \mathbf{R}^{3})}.
\endaligned
\end{equation}
\end{definition}

We will also utilize the local smoothing estimate of \cite{RV}. Suppose $\psi$ is the same $\psi$ as in Definition $\ref{d7.1}$. Then

\begin{equation}\label{7.19}
\| \psi(\frac{x}{R}) e^{it \Delta} (P_{j} u_{0}) \|_{L_{t,x}^{2}(I \times \mathbf{R}^{d})} \lesssim 2^{-j/2} R^{1/2} \| P_{j} u_{0} \|_{L^{2}(\mathbf{R}^{d})}.
\end{equation}

\noindent The dual of $(\ref{7.19})$ is

\begin{equation}\label{7.20}
\| \int_{I} e^{-it \Delta} \psi(\frac{x}{R}) (P_{j} F(\tau)) \|_{L_{x}^{2}(\mathbf{R}^{d})} \lesssim 2^{-j/2} R^{1/2} \| \psi(\frac{x}{R}) P_{j} F \|_{L_{t,x}^{2}(I \times \mathbf{R}^{d})}.
\end{equation}

\noindent Interpolating $(\ref{7.14})$ and $(\ref{7.20})$, for any $q < 2$, if $F$ is supported on $|x| \leq R$,

\begin{equation}\label{3.15}
\| |\nabla|^{1 - \frac{1}{q}} \int_{0}^{\infty} e^{-it \Delta} F(t,x) dt \|_{L^{2}(\mathbf{R}^{3})} \lesssim R^{1 - \frac{1}{q}} \| F \|_{L_{t}^{q} L_{x}^{2}(\mathbf{R} \times \mathbf{R}^{3})}.
\end{equation}

\noindent Now let $\chi(x) = \psi(\frac{x}{2}) - \psi(x)$. For any $0 < R < \infty$ and $x \in \mathbf{R}^{3}$,

\begin{equation}
1 = \psi(R x) + \sum_{j = 0}^{\infty} \chi(2^{-j} Rx).
\end{equation}

\noindent Then by $(\ref{3.15})$,

\begin{equation}\label{3.16}
\| |\nabla|^{1 - \frac{1}{q}} \int_{I} e^{-it \Delta} F(t) dt \|_{L_{x}^{2}(\mathbf{R}^{3})} \lesssim R^{\frac{1}{q} - 1} \| \psi(Rx) F \|_{L_{t}^{q} L_{x}^{2}} + R^{\frac{1}{q} - 1} \sum_{j \geq 0} 2^{j(1 - \frac{1}{q})} \| \chi(2^{-j} Rx) F \|_{L_{t}^{q} L_{x}^{2}}.
\end{equation}

\noindent To simplify notation, let

\begin{equation}\label{3.16.1}
\| F \|_{X_{R}(I \times \mathbf{R}^{d})} = R^{\frac{1}{q} - 1} \| \psi(Rx) F \|_{L_{t}^{q} L_{x}^{2}} + R^{\frac{1}{q} - 1} \sum_{j \geq 0} 2^{j(1 - \frac{1}{q})} \| \chi(2^{-j} Rx) F \|_{L_{t}^{q} L_{x}^{2}}.
\end{equation}

\subsection{$U_{\Delta}^{2}$ spaces}
\noindent The $U_{\Delta}^{2}$ spaces are a class of function spaces first introduced in \cite{Tat} to study wave maps. \cite{KoTa1} and \cite{KoTa2} applied these spaces to nonlinear Schr{\"o}dinger problems. See \cite{HHK} for a general description of these spaces. These spaces are quite useful to critical problems since the $X^{s,b}$ spaces of \cite{B} and \cite{B1} (see also \cite{Gin}) are not scale invariant except at $b = \frac{1}{2}$, which has the same difficulty as the failure of the embedding $\dot{H}^{1/2}(\mathbf{R}) \subset L^{\infty}(\mathbf{R})$.

\begin{definition}[$U_{\Delta}^{p}$ spaces]\label{d7.6}
Let $1 \leq p < \infty$. Let $U_{\Delta}^{p}$ be an atomic space whose atoms are piecewise solutions to the linear equation,

\begin{equation}\label{7.21}
u_{\lambda} = \sum_{k} 1_{[t_{k}, t_{k + 1})} e^{it \Delta} u_{k}, \hspace{5mm} \sum_{k} \| u_{k} \|_{L^{2}}^{p} = 1.
\end{equation}

\noindent Then for any $1 \leq p < \infty$,

\begin{equation}\label{7.22}
\| u \|_{U_{\Delta}^{p}} = \inf \{ \sum_{\lambda} |c_{\lambda}| : u = \sum_{\lambda} c_{\lambda} u_{\lambda}, \text{$u_{\lambda}$ are $U_{\Delta}^{p}$ atoms} \}.
\end{equation}
\end{definition}

\noindent For any $1 \leq p < \infty$, $U_{\Delta}^{p} \subset L^{\infty} L^{2}$. Additionally, $U_{\Delta}^{p}$ functions are continuous except at countably many points and right continuous everywhere.

%\begin{definition}[$V_{\Delta}^{p}$ spaces]\label{d7.7}
%Let $1 \leq p < \infty$. Then $V_{\Delta}^{p}$ is the space of right continuous functions $u \in L^{\infty}(L^{2})$ such that

%\begin{equation}\label{7.23}
%\| v \|_{V_{\Delta}^{p}}^{p} = \| v \|_{L^{\infty}(L^{2})}^{p} + \sup_{\{ t_{k} \} \nearrow} \sum_{k} \| e^{-it_{k} \Delta} v(t_{k}) - e^{-it_{k + 1} \Delta} v(t_{k + 1}) \|_{L^{2}}^{p}.
%\end{equation}

%\noindent The supremum is taken over increasing sequences $t_{k}$.
%\end{definition}

\begin{theorem}\label{t7.8}
If $u$ solves

\begin{equation}\label{7.24}
i u_{t} + \Delta u = F_{1} + F_{2}, \hspace{5mm} u(0, x) = u_{0},
\end{equation}

\noindent on the interval $0 \in I \subset \mathbf{R}$, then for $q < 2$,

\begin{equation}\label{7.25}
\| |\nabla|^{1 - \frac{1}{q}} u \|_{U_{\Delta}^{2}(I \times \mathbf{R}^{d})} \lesssim_{q} \| |\nabla|^{1 - \frac{1}{q}} u_{0} \|_{L^{2}(\mathbf{R}^{d})} + \| F_{1} \|_{X_{R}(I \times \mathbf{R}^{d})} + \| |\nabla|^{1 - \frac{1}{q}} F_{2} \|_{L_{t}^{2+} L_{x}^{\frac{2d}{d - 2}-}(I \times \mathbf{R}^{d})}.
\end{equation}

%\noindent By Duhamel's formula

%\begin{equation}\label{7.26}
%\| u \|_{U_{\Delta}^{p}} \lesssim \| u(0) \|_{L^{2}} + \| (i \partial_{t} + \partial_{x}^{2}) u \|_{DU_{\Delta}^{p}}.
%\end{equation}

%\noindent Finally, there is the duality relation

%\begin{equation}\label{7.27}
%(DU_{\Delta}^{p})^{\ast} = V_{\Delta}^{p'}.
%\end{equation}

%\noindent These spaces are also closed under truncation in time.

%\begin{equation}\label{7.28}
%\aligned
%\chi_{I} : U_{\Delta}^{p} \rightarrow U_{\Delta}^{p}, \\
%\chi_{I} : V_{\Delta}^{p} \rightarrow V_{\Delta}^{p}.
%\endaligned
%\end{equation}
\end{theorem}

\noindent \emph{Proof:} This is proved using Strichartz estimates, $(\ref{7.20})$, and the Christ-Kiselev lemma (see \cite{ChrKis}). $\Box$\medskip

\noindent \textbf{Remark:} The notation
\begin{equation}
\| A \|_{L_{t}^{p+} L_{x}^{q-}} \lesssim B,
\end{equation}
means that for admissible pair $(\tilde{p}, \tilde{q})$ close to $(p, q)$ with $\tilde{p} > p$ and $\tilde{q} < q$,
\begin{equation}
\| A \|_{L_{t}^{\tilde{p}} L_{x}^{\tilde{q}}} \lesssim B,
\end{equation}
for an admissible pair $(\tilde{p}, \tilde{q})$, where the implicit constant can go to infinity as $(\tilde{p}, \tilde{q}) \rightarrow (p, q)$.

\section{Description of the I-method and outline of the proof}
\noindent Since there are no known conserved quantities that control $\| u \|_{\dot{H}^{s}}$ for $0 < s < 1$, we utilize the by now well known modified energy of \cite{CKSTT1}, $E(Iu(t))$, where $I$ is the $I$-operator and $E$ is the usual energy.

\begin{definition}[I-operator]\label{d2.1}
\noindent Let $I : H^{s}(\mathbf{R}^{d}) \rightarrow H^{1}(\mathbf{R}^{d})$ be the Fourier multiplier

\begin{equation}\label{2.1}
\widehat{If}(\xi) = m_{N}(\xi) \hat{f}(\xi),
\end{equation}

\noindent where

\begin{equation}\label{2.2}
m_{N}(\xi) =
\left\{
	\begin{array}{ll}
		1  & \mbox{if } |\xi| \leq N, \\
		\frac{N^{1 - s}}{|\xi|^{1 - s}} & \mbox{if } |\xi| \geq 2N.
	\end{array}
\right.
\end{equation}
\end{definition}

\noindent To simplify notation $N$ is suppressed for the rest of the paper.\vspace{5mm}

\noindent \textbf{Remark:} It is convenient to write

\begin{equation}
P_{\leq N} u = \sum_{j : 2^{j} \leq N} P_{j} u.
\end{equation}

\noindent Let $P_{> N} u = u - P_{\leq N} u$. This notation may be abbreviated $u_{\leq N} = P_{\leq N} u$.\vspace{5mm}

\noindent There is an obvious tradeoff here. Taking $N = \infty$, $E(Iu(t))$ is the energy of $u$, which is a conserved quantity. However, for a general $u \in \dot{H}^{s}$, $s < 1$, $E(Iu(0)) = \infty$. In general, as $N$ increases, $\frac{d}{dt} E(Iu(t))$ decreases and $E(Iu(t))$ increases. Therefore, the question of global well-posedness revolves around which side will win this tug of war. More precisely, by the Sobolev embedding theorem, when $d = 3$,

\begin{equation}\label{2.3}
E(Iu(t)) \lesssim \| Iu \|_{\dot{H}^{1}(\mathbf{R}^{3})}^{2} + \| Iu \|_{\dot{H}^{1}(\mathbf{R}^{3})}^{2} \| u \|_{\dot{H}^{1/2}(\mathbf{R}^{3})}^{2},
\end{equation}

\noindent and when $d = 2$,

\begin{equation}\label{2.4}
E(Iu(t)) \lesssim \| Iu \|_{\dot{H}^{1}(\mathbf{R}^{2})}^{2} + \| Iu \|_{\dot{H}^{1}(\mathbf{R}^{2})}^{2} \| u \|_{\dot{H}^{1 - \frac{1}{k}}(\mathbf{R}^{2})}^{2k}.
\end{equation}

\noindent Therefore,

\begin{equation}\label{2.5}
E(Iu(0)) \lesssim C(\| u(0) \|_{H^{s}}) N^{2(1 - s)}.
\end{equation}

\noindent Meanwhile,

\begin{equation}\label{2.6}
 \| u(t) \|_{H^{s}(\mathbf{R}^{d})}^{2} \lesssim E(Iu(t)) + M(Iu(t)).
\end{equation}

\noindent Since $M(Iu(t)) \leq M(u(t)) = M(u(0))$, a uniform bound on $E(Iu(t))$ for all $t$ yields a uniform bound on $\| u(t) \|_{H^{s}(\mathbf{R}^{3})}$.\vspace{5mm}

\noindent It is convenient to use the rescaling in $(\ref{1.3})$ so that $E(Iu(t)) \leq 1$. Indeed, by $(\ref{1.4})$, when $d = 2$ there exists $\lambda^{s - 1 + \frac{1}{k}} \sim C(\| u(0) \|_{H^{s}(\mathbf{R}^{2})}) N^{s - 1}$ and when $d = 3$ there exists $\lambda^{s - \frac{1}{2}} \sim C(\| u(0) \|_{H^{s}(\mathbf{R}^{3})}) N^{s - 1}$ such that

\begin{equation}\label{2.7}
 E(Iu_{\lambda}(0)) \leq \frac{1}{2}.
\end{equation}

\noindent \textbf{Remark:} $C(\| u(0) \|_{H^{s}(\mathbf{R}^{d})})$ is a constant that may change from line to line.\vspace{5mm}

\noindent Then by $(\ref{1.4})$,

\begin{equation}\label{2.8}
\| u_{\lambda}(0) \|_{L_{x}^{2}(\mathbf{R}^{3})} \lesssim C(\| u(0) \|_{H^{s}(\mathbf{R}^{3})}) N^{\frac{1 - s}{2s - 1}} \| u(0) \|_{L_{x}^{2}(\mathbf{R}^{3})},
\end{equation}

\noindent and

\begin{equation}\label{2.8.1}
\| u_{\lambda}(0) \|_{L_{x}^{2}(\mathbf{R}^{2})} \lesssim C(\| u(0) \|_{H^{s}(\mathbf{R}^{2})}) N^{\frac{k - 1}{k} \cdot \frac{1 - s}{s - 1 + \frac{1}{k}}} \| u(0) \|_{L_{x}^{2}(\mathbf{R}^{2})}.
\end{equation}

\noindent $\lambda$ is suppressed until the end of the paper, so for now $u$ refers to $u_{\lambda}$ until otherwise indicated.\vspace{5mm}

Next recall the interaction Morawetz estimate.

\begin{theorem}[Interaction Morawetz estimate]\label{t2.2}
Suppose $u$ is a solution to $(\ref{1.1})$ or $(\ref{1.2})$ on some interval $J$. Then

\begin{equation}\label{2.9}
\| |\nabla|^{\frac{3 - d}{2}} |u|^{2} \|_{L_{t,x}^{2}(J \times \mathbf{R}^{d})}^{2} \lesssim \| u \|_{L_{t}^{\infty} L_{x}^{2}(J \times \mathbf{R}^{d})}^{2} \| u \|_{L_{t}^{\infty} \dot{H}_{x}^{1/2}(J \times \mathbf{R}^{d})}^{2}.
\end{equation}
\end{theorem}

\noindent \emph{Proof:} This was proved in three dimensions by \cite{CKSTT2}. \cite{CGT1} and \cite{PV} independently proved $(\ref{2.9})$ in dimensions one and two. \cite{TVZ} proved the interaction Morawetz estimate in dimensions $d \geq 4$, a result that will not be needed here. $\Box$\vspace{5mm}

\noindent $(\ref{2.9})$ is extremely useful due to a local well-posedness result of \cite{CKSTT2}.

\begin{lemma}\label{l2.3}
 If $E(Iu(a_{l})) \leq 1$, $J_{l} = [a_{l}, b_{l}]$, and $\| u \|_{L_{t,x}^{4}(J_{l} \times \mathbf{R}^{3})} \leq \epsilon$ for some $\epsilon > 0$ sufficiently small, then

\begin{equation}\label{2.9.1}
\| \nabla Iu \|_{S^{0}(J_{l} \times \mathbf{R}^{3})} \lesssim 1.
\end{equation}
\end{lemma}

\noindent \emph{Proof:} See \cite{CKSTT2} or \cite{D1}. $\Box$\vspace{5mm}

\noindent A similar result is available in dimension $d = 2$.

\begin{lemma}\label{l2.3}
 If $k > 1$, $E(Iu(a_{l})) \leq 1$, $J_{l} = [a_{l}, b_{l}]$, and $\| |\nabla|^{1/2} |u|^{2} \|_{L_{t,x}^{2}(J_{l} \times \mathbf{R}^{2})} \leq \epsilon$ for some $\epsilon(k) > 0$ sufficiently small, then for $(p, q) \in \mathcal A_{2}$,

\begin{equation}\label{2.10}
\| \nabla Iu \|_{L_{t}^{p} L_{x}^{q}(J_{l} \times \mathbf{R}^{2})} \lesssim_{p, k} 1.
\end{equation}
\end{lemma}

\noindent \emph{Proof:} By the Sobolev embedding theorem

\begin{equation}\label{2.11}
\| u \|_{L_{t}^{4} L_{x}^{8}(J_{l} \times \mathbf{R}^{2})}^{2} \lesssim \| |\nabla|^{1/2} |u|^{2} \|_{L_{t,x}^{2}(J_{l} \times \mathbf{R}^{2})} \leq \epsilon.
\end{equation}

\noindent Interpolating $(\ref{2.11})$ with $\| P_{j} u \|_{L_{x}^{\infty}(J_{l} \times \mathbf{R}^{2})} \lesssim 2^{j} \| P_{j} u \|_{L_{x}^{2}(\mathbf{R}^{2})}$, combined with the Littlewood-Paley theorem proves

\begin{equation}\label{2.12}
\| Iu \|_{L_{t}^{3k} L_{x}^{6k}(J_{l} \times \mathbf{R}^{2})} \lesssim \epsilon^{\frac{4}{3k}} \| \nabla Iu \|_{L_{t}^{\infty} L_{x}^{2}(J_{l} \times \mathbf{R}^{2})}^{1 - \frac{4}{3k}}.
\end{equation}

\noindent Also by Bernstein's inequality and $(\ref{2.2})$,

\begin{equation}\label{2.13}
\| (1 - I)u \|_{L_{t}^{3k} L_{x}^{6k}(J_{l} \times \mathbf{R}^{2})} \lesssim N^{-\frac{1}{k}} \| \nabla Iu \|_{L_{t}^{3k} L_{x}^{\frac{6k}{3k - 2}}(J_{l} \times \mathbf{R}^{2})}.
\end{equation}

\noindent Then by Strichartz estimates (Theorem $\ref{t7.4}$),

\begin{equation}\label{2.14}
\aligned
\| \nabla Iu \|_{L_{t}^{3k} L_{x}^{\frac{6k}{3k - 2}} \cap L_{t}^{\infty} L_{x}^{2}(J_{l} \times \mathbf{R}^{2})} \lesssim \| \nabla Iu(a_{l}) \|_{L_{x}^{2}(\mathbf{R}^{2})}   \\ + (\epsilon^{\frac{4}{3k}} \| \nabla Iu \|_{L_{t}^{\infty} L_{x}^{2}(J_{l} \times \mathbf{R}^{2})}^{\frac{3k - 4}{3k}} + N^{-\frac{1}{k}} \| \nabla Iu \|_{L_{t}^{3k} L_{x}^{\frac{6k}{3k - 2}}(J_{l} \times \mathbf{R}^{2})})^{2k} \| \nabla Iu \|_{L_{t}^{\infty} L_{x}^{2}(J_{l} \times \mathbf{R}^{2})}.
\endaligned
\end{equation}

\noindent Since $N$ is large and $\epsilon > 0$ is small the proof is complete. $\Box$\vspace{5mm}

\noindent $(\ref{2.14})$ also implies

\begin{equation}\label{2.15}
\| \nabla Iu \|_{U_{\Delta}^{2}(J_{l} \times \mathbf{R}^{2})} \lesssim_{k} 1,
\end{equation}

\noindent and similarly Lemma $\ref{l2.3}$ implies

\begin{equation}\label{2.16}
\| \nabla Iu \|_{U_{\Delta}^{2}(J_{l} \times \mathbf{R}^{3})} \lesssim 1.
\end{equation}

\noindent Theorems $\ref{t1.4}$, $\ref{t1.4.1}$, and $\ref{t1.5}$ are then proved by a bootstrapping estimate. Let

\begin{equation}\label{2.17}
J = \{ t : E(Iu(\tau)) \leq 1 \hspace{5mm} \text{for all} \hspace{5mm} 0 \leq \tau \leq t \}.
\end{equation}

\noindent $J$ is clearly nonempty since $0 \in J$. Moreover, standard local well-posedness theory implies that $J$ is a closed interval. Therefore, to prove $J = [0, \infty)$ it suffices to show that $J$ is open. By $(\ref{2.2})$, interpolation, and Bernstein's inequality,

\begin{equation}\label{2.18}
\| P_{\leq N} u(t) \|_{\dot{H}^{1/2}(\mathbf{R}^{d})} \lesssim \| Iu(t) \|_{\dot{H}^{1}(\mathbf{R}^{d})}^{1/2} \| P_{\leq N} u(t) \|_{L^{2}(\mathbf{R}^{d})}^{1/2},
\end{equation}

\noindent and

\begin{equation}\label{2.19}
\| P_{> N} u(t) \|_{\dot{H}^{1/2}(\mathbf{R}^{d})} \lesssim N^{-1/2} \| Iu \|_{\dot{H}^{1}(\mathbf{R}^{d})}.
\end{equation}

\noindent Therefore if $J$ is an interval such that $E(Iu(t)) \leq 1$ on $J$, then $(\ref{2.8})$, $(\ref{2.9})$, $(\ref{2.18})$, $(\ref{2.19})$, and the conservation of mass imply that

\begin{equation}\label{2.20}
\| u \|_{L_{t,x}^{4}(J \times \mathbf{R}^{3})}^{4} \lesssim C(\| u(0) \|_{H^{s}(\mathbf{R}^{3})}) N^{\frac{3(1 - s)}{2s - 1}},
\end{equation}

\noindent and

\begin{equation}\label{2.20.1}
\| u \|_{L_{t}^{4} L_{x}^{8}(J \times \mathbf{R}^{2})}^{4} \lesssim C(\| u(0) \|_{H^{s}(\mathbf{R}^{2})}, k) N^{\frac{k - 1}{k} \cdot \frac{3(1 - s)}{s - 1 + \frac{1}{k}}}. 
\end{equation}

\noindent To close the bootstrap, we prove that for $N(d, k, \| u(0) \|_{H^{s}})$ sufficiently large,

\begin{equation}\label{2.23}
\int_{J} \frac{d}{dt} E(Iu(t)) dt \leq \frac{1}{10}.
\end{equation}

\noindent $(\ref{2.23})$ and Theorem $\ref{t1.1}$ imply that for any $T > 0$ there exists $\delta(T) > 0$ such that if $[0, T] \subset J$, $[0, T + \delta) \subset J$. Therefore $J$ is open and thus $J = [0, \infty)$. Finally, we can recover the $\| u(t) \|_{H^{s}}$ bound by rescaling back and then computing the $\| u(t) \|_{H^{s}}$ norm from the bounds on $M(u(t))$ and $E(Iu(t))$ after rescaling.\vspace{5mm}

\noindent $(\ref{2.23})$ is proved using long time Strichartz estimates. Estimates of this form were introduced in \cite{D2} within the context of the mass-critical nonlinear Schr{\"o}dinger initial value problem. The long time Strichartz estimates have been utilized in subsequent papers (\cite{D3}, \cite{D4}, \cite{KV2}, \cite{Murphy2}, \cite{Murphy1}, \cite{Murphy3}, \cite{Visan1}).

%\noindent Long time  When $d = 3$ partition $J$ into $\sim_{\| u(0) \|_{H^{s}}} N^{\frac{3(1 - s)}{2s - 1}}$ subintervals such that on each subinterval $J_{l}$,

%\begin{equation}\label{2.21}
% \| u \|_{L_{t,x}^{4}(J_{l} \times \mathbf{R}^{3})} \leq \epsilon.
%\end{equation}

%\noindent When $d = 2$ partition $J$ into $\sim_{\| u(0) \|_{H^{s}}, k} N^{\frac{k - 1}{k} \cdot \frac{3(1 - s)}{s - 1 + \frac{1}{k}}}$ subintervals such that

%\begin{equation}\label{2.22}
%\| u \|_{L_{t}^{4} L_{x}^{8}(J_{l} \times \mathbf{R}^{2})} \leq \epsilon.
%\end{equation}

%\noindent 

\section{Induction on frequency and long time Strichartz estimates in three dimensions}
%Examining $(\ref{2.2})$, it is convenient to define the notation

%\begin{equation}\label{3.0}
%P_{> N} u = \sum_{2^{j} > N} P_{j} u, \hspace{5mm} P_{\leq N} u = \sum_{2^{j} \leq N} P_{j} u,
%\end{equation}

%\noindent as well as

%\begin{equation}\label{3.0.1}
%P_{> M} u = \sum_{2^{j} > M} P_{j} u, \hspace{5mm} P_{\leq M} u = \sum_{2^{j} \leq M} P_{j} u.
%\end{equation}

%\noindent We prove

\begin{theorem}\label{t3.1}
Let $0 \in J$ be an interval such that $E(Iu(t)) \leq 1$ on $J$. Then for $N(s, \| u_{0} \|_{H^{s}})$ sufficiently large,

\begin{equation}\label{3.0.2}
\| P_{> \frac{N}{8}} \nabla Iu \|_{L_{t}^{2} L_{x}^{6}(J \times \mathbf{R}^{3})} \lesssim 1.
\end{equation}
\end{theorem}

\noindent \emph{Proof:} As in \cite{D2}, this theorem is proved using an induction on frequency argument. First observe that

\begin{equation}\label{3.3}
P_{> M} (|u_{\leq \frac{M}{8}}|^{2} u_{\leq \frac{M}{8}}) = 0.
\end{equation}

\noindent \textbf{Remark:} This fact is why this method does not immediately carry over to a non-algebraic nonlinearity, in other words, when $p$ is not equal to $2k$ for some positive integer $k$.\vspace{5mm}

\noindent Decompose

\begin{equation}\label{3.3.1}
P_{> M} F(u) = P_{> M} O((u_{> \frac{M}{8}}) (u_{\leq \frac{M}{8}})^{2}) + P_{> M} O((u_{> \frac{M}{8}})^{2} u).
\end{equation}

\noindent By the product rule and the fact that $\nabla I$ is a Fourier multiplier whose symbol is increasing as $|\xi| \nearrow \infty$, if $M \leq N$,

\begin{equation}\label{3.3.2}
\nabla I P_{> M} O((u_{> \frac{M}{8}})(u_{\leq \frac{M}{8}})^{2}) = O((\nabla I P_{> \frac{M}{8}} u)(P_{\leq \frac{M}{8}} u)^{2}) + O((I P_{> \frac{M}{8}} u)(\nabla u_{\leq \frac{M}{8}})(u_{\leq \frac{M}{8}})).
\end{equation}

%\noindent Now take the cutoff supported on the annulus $|x| \sim 2^{j}$, $\psi_{j}(x) = \chi^{2}(2^{-j} x) - \chi^{2}(2^{-j + 1} x)$, where $j \geq 0$.

\noindent Then by $(\ref{7.25})$,

\begin{equation}\label{3.1}
\aligned
\| \nabla I P_{> M} u(t) \|_{U_{\Delta}^{2}(J \times \mathbf{R}^{3})} \lesssim \| \nabla I P_{> M} u(0) \|_{L_{x}^{2}(\mathbf{R}^{3})} + \| \nabla I P_{> M} O((u_{> \frac{M}{8}})^{2} u) \|_{L_{t}^{2-} L_{x}^{6/5+} (J \times \mathbf{R}^{3})} \\
+ \| P_{> M} O((u_{> \frac{M}{8}})(\nabla u_{\leq \frac{M}{8}})(u_{\leq \frac{M}{8}})) \|_{L_{t}^{2-} L_{x}^{6/5+}} + \frac{1}{M^{1 - \frac{1}{q}}} \| P_{> M} O((\nabla I u_{> \frac{M}{8}})(u_{\leq \frac{M}{8}})^{2}) \|_{X_{R}},
\endaligned
\end{equation}

\noindent for some $R$ to be specified later.\medskip

First observe that since $E(Iu(t)) \leq 1$ for all $t \in J$,

\begin{equation}
\| \nabla I P_{> M} u(0) \|_{L_{x}^{2}(\mathbf{R}^{3})} \lesssim 1.
\end{equation}

\noindent Again using the properties of $\nabla I$, choosing $\delta(\epsilon) > 0$ so that $(2 - \epsilon, \frac{6}{5} + \delta(\epsilon))$ is the dual of an admissible pair, and subsequent $\epsilon$ and $\delta(\epsilon)$ to correspond with H\"older's inequality,

\begin{equation}\label{3.4}
\| \nabla I((P_{> \frac{M}{8}} u)^{2} u) \|_{L_{t}^{2 - \epsilon} L_{x}^{6/5 + \delta(\epsilon)}(J \times \mathbf{R}^{3})} \lesssim \| \nabla Iu \|_{L_{t}^{\infty - \epsilon} L_{x}^{2 + \delta(\epsilon)}(J \times \mathbf{R}^{3})} \| P_{> \frac{M}{8}} u \|_{L_{t}^{4} L_{x}^{6}(J \times \mathbf{R}^{3})}^{2}
\end{equation}

\begin{equation}\label{3.5}
+ \| \nabla I P_{> \frac{M}{8}} u \|_{L_{t}^{2} L_{x}^{6}(J \times \mathbf{R}^{3})} \| P_{> \frac{M}{8}} u \|_{L_{t}^{\infty} L_{x}^{2}(J \times \mathbf{R}^{3})} \| P_{\leq N} u \|_{L_{t}^{\infty - \epsilon} L_{x}^{6 + \delta(\epsilon)}(J \times \mathbf{R}^{3})}
\end{equation}

\begin{equation}\label{3.6}
+ \| \nabla I P_{> \frac{M}{8}} u \|_{L_{t}^{2} L_{x}^{6}(J \times \mathbf{R}^{3})} \| P_{> \frac{M}{8}} u \|_{L_{t}^{\infty} L_{x}^{3}(J \times \mathbf{R}^{3})} \| P_{> N} u \|_{L_{t}^{\infty - \epsilon} L_{x}^{3 + \delta(\epsilon)}(J \times \mathbf{R}^{3})}.
\end{equation}
\noindent \textbf{Remark:} The notation $\infty - \epsilon$ refers to a very large number, specifically $\frac{2(2 - \epsilon)}{\epsilon}$.\medskip

%\begin{equation}\label{3.7}
%+ \| \nabla Iu \|_{L_{t}^{\infty-} L_{x}^{2+}(J \times \mathbf{R}^{3})} \| P_{> \frac{M}{8}} u \|_{L_{t}^{4} L_{x}^{6}(J \times \mathbf{R}^{3})}^{2}.
%\end{equation}

Now by the Sobolev embedding theorem and interpolation,

\begin{equation}\label{3.10}
\| P_{> \frac{M}{8}} u \|_{L_{t}^{4} L_{x}^{6}(J \times \mathbf{R}^{3})} \lesssim \| |\nabla|^{1/2} P_{> \frac{M}{8}} u \|_{L_{t}^{2} L_{x}^{6}(J \times \mathbf{R}^{3})}^{1/2} \| |\nabla|^{1/2} P_{> \frac{M}{8}} u \|_{L_{t}^{\infty} L_{x}^{2}(J \times \mathbf{R}^{3})}^{1/2}.
\end{equation}

\noindent Therefore, by Bernstein's inequality, $E(Iu(t)) \leq 1$, and $(\ref{2.2})$,

\begin{equation}\label{3.11}
(\ref{3.10}) \lesssim M^{-1/2} \| \nabla I P_{> \frac{M}{8}} u \|_{L_{t}^{2} L_{x}^{6}(J \times \mathbf{R}^{3})}^{1/2}.
\end{equation}

\noindent Therefore,

\begin{equation}\label{3.11.1}
\aligned
(\ref{3.4}) + (\ref{3.5}) + (\ref{3.6}) \lesssim \| \nabla I P_{> \frac{M}{8}} u \|_{L_{t}^{2} L_{x}^{6}(J \times \mathbf{R}^{3})} \\ \times (M^{-1} \| \nabla Iu \|_{L_{t}^{\infty - \epsilon} L_{x}^{2 + \delta(\epsilon)}} + M^{-1} \| P_{\leq N} u \|_{L_{t}^{\infty - \epsilon} L_{x}^{6 + \delta(\epsilon)}} + M^{-1/2} \| P_{> N} u \|_{L_{t}^{\infty - \epsilon} L_{x}^{3 + \delta(\epsilon)}}).
\endaligned
\end{equation}

\noindent Now by $(\ref{2.16})$, $(\ref{2.20})$, $E(Iu(t)) \leq 1$ on $J$, and Lemma $\ref{l2.3}$,

\begin{equation}\label{3.8}
\| Iu \|_{L_{t}^{\infty - \epsilon} L_{x}^{6 + \delta(\epsilon)}(J \times \mathbf{R}^{3})} + \| \nabla Iu \|_{L_{t}^{\infty - \epsilon} L_{x}^{2 + \delta(\epsilon)}(J \times \mathbf{R}^{3})}  \lesssim_{s, \| u_{0} \|_{H^{s}}(\mathbf{R}^{3})} N^{\frac{3(1 - s)}{2s - 1} \cdot \frac{\epsilon}{2(2 - \epsilon)}}
\end{equation}

\noindent and

\begin{equation}\label{3.9}
\| P_{> N} u \|_{L_{t}^{\infty - \epsilon} L_{x}^{3 + \delta(\epsilon)}(J \times \mathbf{R}^{3})} \lesssim N^{-\frac{1}{2} + \frac{3(1 - s)}{2s - 1} \cdot \frac{\epsilon}{2(2 - \epsilon)}}.
\end{equation}

\noindent Therefore, if $M \leq N$,

\begin{equation}\label{3.12}
(\ref{3.11.1}) \lesssim_{s, \| u_{0} \|_{H^{s}(\mathbf{R}^{3})}} M^{-1} N^{\frac{3(1 - s)}{2s - 1} \cdot \frac{\epsilon}{2(2 - \epsilon)}} \| \nabla I P_{> \frac{M}{8}} u \|_{U_{\Delta}^{2}(J \times \mathbf{R}^{3})}.
\end{equation}

\noindent Similarly,

\begin{equation}\label{3.12.1}
\aligned
\| P_{> M} O((u_{> \frac{M}{8}})(\nabla u_{\leq \frac{M}{8}})(u_{\leq \frac{M}{8}})) \|_{L_{t}^{2 - \epsilon} L_{x}^{6/5 + \delta(\epsilon)}} \\ \lesssim \|  \nabla u_{\leq \frac{M}{8}} \|_{L_{t}^{\infty - \epsilon} L_{x}^{2 + \delta(\epsilon)}} \| u_{> \frac{M}{8}} \|_{L_{t}^{2} L_{x}^{6}} \| u_{\leq \frac{M}{8}} \|_{L_{t}^{\infty} L_{x}^{6}} \\
\lesssim_{s, \| u_{0} \|_{H^{s}(\mathbf{R}^{3})}} M^{-1} N^{\frac{3(1 - s)}{2s - 1} \cdot \frac{\epsilon}{2(2 - \epsilon)}} \| \nabla I P_{> \frac{M}{8}} u \|_{U_{\Delta}^{2}(J \times \mathbf{R}^{3})}.
\endaligned
\end{equation}

It only remains to analyze

\begin{equation}\label{3.13}
\frac{1}{M^{1 - \frac{1}{q}}} \| (\nabla I P_{> \frac{M}{8}} u)(u_{\leq \frac{M}{8}})^{2} \|_{X_{R}}.
\end{equation}
\textbf{Remark:} For notational convenience, choose $q = 2 - \epsilon$.\medskip

\noindent It is here that the radial symmetry of $u$ is utilized. Recall that for any $\frac{1}{2} < s < \frac{3}{2}$, the radial Sobolev embedding implies

\begin{equation}\label{3.18}
\aligned
\| |x|^{3/2 - s} u \|_{L_{x}^{\infty}(\mathbf{R}^{3})} \lesssim  \| u \|_{\dot{H}^{s}(\mathbf{R}^{3})}.
\endaligned
\end{equation}

\noindent Interpolating this with

\begin{equation}\label{3.19}
\| Iu \|_{L_{t}^{4} L_{x}^{\infty}(J \times \mathbf{R}^{3})}^{4} \lesssim_{\| u_{0} \|_{H^{s}}(\mathbf{R}^{3})} N^{\frac{3(1 - s)}{2s - 1}},
\end{equation}

\noindent which is a consequence of $(\ref{2.16})$ and Strichartz estimates, along with $(\ref{2.8})$, implies that

\begin{equation}\label{3.20}
\| |x|^{1/2} Iu \|_{L_{t}^{\infty - \epsilon} L_{x}^{\infty}(J \times \mathbf{R}^{3})} \lesssim_{s, \| u_{0} \|_{H^{s}}(\mathbf{R}^{3})} N^{\frac{3(1 - s)}{2s - 1} \cdot \frac{\epsilon}{2(2 - \epsilon)}} N^{\frac{\epsilon}{2 - 3 \epsilon} \cdot \frac{1 - s}{2s - 1}}.
\end{equation}

\noindent Now choose $R = N$. By $(\ref{7.19})$, $(\ref{3.20})$, and $(\ref{3.18})$,

\begin{equation}\label{3.21}
\aligned
R^{\frac{1}{q} - 1} M^{\frac{1}{q} - 1} \| \psi(Rx) (\nabla I P_{> \frac{M}{8}} u)(u_{\leq \frac{M}{8}})^{2} \|_{L_{t}^{q} L_{x}^{2}(J \times \mathbf{R}^{3})} \\ 
\lesssim R^{\frac{1}{q} - 1} M^{\frac{1}{q} - 1} \| \psi(Rx) (\nabla I P_{> \frac{M}{8}} u) \|_{L_{t,x}^{2}(J \times \mathbf{R}^{3})} \| u_{\leq \frac{M}{8}} \|_{L_{t}^{4} L_{x}^{\infty}}^{\frac{2 \epsilon}{(2 - \epsilon)}} \| u_{\leq \frac{M}{8}} \|_{L_{t,x}^{\infty}}^{\frac{4 - 4 \epsilon}{2 - \epsilon}} \\
\lesssim_{s, \| u_{0} \|_{H^{s}}} R^{\frac{1}{q} - 1} R^{-\frac{1}{2}} M^{\frac{1}{q} - 1} M^{-\frac{1}{2}} N^{\frac{3(1 - s)}{2s - 1} \cdot \frac{2 \epsilon}{2 - \epsilon}} M^{\frac{2 - 2 \epsilon}{2 - \epsilon}} \| \nabla I P_{> \frac{M}{8}} u \|_{U_{\Delta}^{2}(J \times \mathbf{R}^{3})} \\
= N^{\frac{-4 + 3 \epsilon}{2(2 - \epsilon)}} M^{\frac{-4 + 3 \epsilon}{2(2 - \epsilon)}} N^{\frac{3(1 - s)}{2s - 1} \cdot \frac{2 \epsilon}{2 - \epsilon}} M^{\frac{2 - 2 \epsilon}{2 - \epsilon}}  \| \nabla I P_{> \frac{M}{8}} u \|_{U_{\Delta}^{2}(J \times \mathbf{R}^{3})}.
\endaligned
\end{equation}

\noindent Also, by $(\ref{7.19})$ and $(\ref{3.20})$,

\begin{equation}\label{3.22}
\aligned
M^{\frac{1}{q} - 1} \sum_{j \geq 0} R^{\frac{1}{q} - 1} 2^{j(1 - \frac{1}{q})} \| \chi(2^{-j} Rx) (\nabla I P_{> \frac{M}{8}} u)(u_{\leq \frac{M}{8}})^{2} \|_{L_{t}^{q} L_{x}^{2}} \\ 
\lesssim M^{\frac{1}{q} - 1} \sum_{j \geq 0} R 2^{-j} R^{\frac{1}{q} - 1} 2^{j(1 - \frac{1}{q})} \| \chi(2^{-j} Rx) (\nabla I P_{> \frac{M}{8}} u) \|_{L_{t,x}^{2}} \| |x|^{1/2} Iu \|_{L_{t}^{\infty - \epsilon} L_{x}^{\infty}} \| |x|^{1/2} Iu \|_{L_{t,x}^{\infty}} \\
\lesssim_{s, \| u_{0} \|_{H^{s}}(\mathbf{R}^{3})} R^{\frac{1}{q} - \frac{1}{2}} M^{\frac{1}{q} - \frac{3}{2}} \sum_{j \geq 0} 2^{j(\frac{1}{2} - \frac{1}{q})} N^{\frac{3(1 - s)}{2s - 1} \cdot \frac{\epsilon}{2(2 - \epsilon)}} N^{\frac{\epsilon}{2 - 3 \epsilon} \cdot \frac{1 - s}{2s - 1}}  \| \nabla I P_{> \frac{M}{8}} u \|_{U_{\Delta}^{2}(J \times \mathbf{R}^{3})} \\ 
\lesssim N^{\frac{\epsilon}{2(2 - \epsilon)}} M^{\frac{\epsilon}{2(2 - \epsilon)} - 1} N^{\frac{3(1 - s)}{2s - 1} \cdot \frac{\epsilon}{2(2 - \epsilon)}} N^{\frac{\epsilon}{2 - 3 \epsilon} \cdot \frac{1 - s}{2s - 1}} \| \nabla I P_{> \frac{M}{8}} u \|_{U_{\Delta}^{2}(J \times \mathbf{R}^{3})}.
\endaligned
\end{equation}

\noindent Combining $(\ref{3.1})$, $(\ref{3.12})$, $(\ref{3.12.1})$, $(\ref{3.21})$, and $(\ref{3.22})$,

\begin{equation}\label{3.23}
\| \nabla I P_{> M} u \|_{U_{\Delta}^{2}(J \times \mathbf{R}^{3})} \lesssim_{s, \| u_{0} \|_{H^{s}(\mathbf{R}^{3})}, \epsilon} 1 + \frac{N^{C_{1}(s) \epsilon}}{M^{1 - C_{2}(s) \epsilon}} \| \nabla I P_{> \frac{M}{8}} u \|_{U_{\Delta}^{2}(J \times \mathbf{R}^{3})}.
\end{equation}
\textbf{Remark:} It is not too important to compute exactly what $C_{1}(s)$ and $C_{2}(s)$ are, except to know that they are constant. This means that for any $s$, after taking $\epsilon(s) > 0$ sufficiently small, $C_{1}(s) \epsilon < \frac{1}{4}$ and $C_{2}(s) \epsilon < \frac{1}{4}$.\medskip

Now we argue by induction on frequency. If $M \geq C(s, \| u_{0} \|_{H^{s}}) N^{2/3}$, then $(\ref{3.23})$ implies
\begin{equation}\label{3.23.1}
\| \nabla I P_{> M} u \|_{U_{\Delta}^{2}(J \times \mathbf{R}^{3})} \lesssim_{s, \| u_{0} \|_{H^{s}(\mathbf{R}^{3})}} 1 + N^{-\frac{1}{4}} C(s, \| u_{0} \|_{H^{s}})^{-3/4} \| \nabla I P_{> \frac{M}{8}} u \|_{U_{\Delta}^{2}(J \times \mathbf{R}^{3})}.
\end{equation}
Also, by $(\ref{2.16})$,
\begin{equation}\label{3.23.2}
\| \nabla I P_{> C(s, \| u_{0} \|_{H^{s}}) N^{2/3}} u \|_{U_{\Delta}^{2}(J \times \mathbf{R}^{3})} \lesssim N^{\frac{3(1 - s)}{4s - 2}}.
\end{equation}
Therefore, by induction, for $C(s, \| u_{0} \|_{H^{s}})$ sufficiently large,
\begin{equation}\label{3.24}
 \| \nabla I P_{> \frac{N}{8}} u \|_{U_{\Delta}^{2}(J \times \mathbf{R}^{3})} \lesssim_{\| u_{0} \|_{H^{s}}, s} 1 + N^{\frac{3(1 - s)}{4s - 2}} N^{-c \ln(N)},
\end{equation}
Therefore, choosing $N$ sufficiently large, say $\ln(N) = C_{0} \frac{1 - s}{s - \frac{1}{2}} + \ln(C(s, \| u_{0} \|_{H^{s}}))$, for some constants $C_{0}$ and $C(s, \| u_{0} \|_{H^{s}})$,

\begin{equation}\label{3.25}
 \| \nabla I P_{> \frac{N}{8}} u \|_{U_{\Delta}^{2}(J \times \mathbf{R}^{3})} \lesssim_{\| u_{0} \|_{H^{s}}(\mathbf{R}^{3})} 1.
\end{equation}

\noindent $\Box$

\section{Energy Increment in three dimensions}
\noindent Now we show a bound on the modified energy increment when $d = 3$.

\begin{lemma}\label{l4.1}
For $N$ sufficiently large so that $\ln(N) \geq C_{0} \frac{1 - s}{s - 1/2} + \ln(C(s, \| u_{0} \|_{H^{s}}))$,

\begin{equation}\label{4.1}
\int_{J} |\frac{d}{dt} E(Iu(t))| dt \lesssim \frac{1}{N^{1-}}.
\end{equation}
\end{lemma}

\noindent \emph{Proof:} Because $I$ is a Fourier multiplier which is constant in time and $\Delta$ commutes with $I$, $(\ref{1.1})$ implies

\begin{equation}\label{4.2}
i Iu_{t} + \Delta Iu = |Iu|^{2} (Iu) + I(|u|^{2} u) - |Iu|^{2} (Iu).
\end{equation}

\noindent Therefore,

\begin{equation}\label{4.3}
\frac{d}{dt} E(Iu(t)) = \langle Iu_{t}, |Iu|^{2}(Iu) - I(|u|^{2} u) \rangle.
\end{equation}

\noindent Then by $(\ref{4.2})$ and integrating by parts,

\begin{equation}\label{4.4}
\aligned
\frac{d}{dt} E(Iu(t)) = -\langle i \nabla Iu, \nabla (|Iu|^{2} (Iu) - I(|u|^{2} u)) \rangle \\
- \langle i I(|u|^{2} u), (|Iu|^{2} (Iu) - I(|u|^{2} u)) \rangle.
\endaligned
\end{equation}

\noindent First estimate $\langle i \nabla Iu, \nabla (|Iu|^{2} (Iu) - I(|u|^{2} u)) \rangle$. As was mentioned before, it is unnecessary to distinguish between polynomial terms involving $u$ and $\bar{u}$. Observe that

\begin{equation}\label{4.6}
 (I P_{\leq \frac{N}{8}} u)^{3} - I((P_{\leq \frac{N}{8}} u)^{3}) = 0.
\end{equation}

\noindent Next,

\begin{equation}\label{4.7}
\aligned
(I P_{> \frac{N}{8}} u)(I P_{\leq \frac{N}{8}} u)^{2} - I((P_{> \frac{N}{8}} u)(P_{\leq \frac{N}{8}} u)^{2}) \\ = (I P_{> \frac{N}{2}} u)(P_{\leq \frac{N}{8}} u)^{2} - I((P_{> \frac{N}{2}} u)(P_{\leq \frac{N}{8}} u)^{2}).
\endaligned
\end{equation}

\noindent By the fundamental theorem of calculus,

\begin{equation}\label{4.8}
 |m(\xi_{2} + \xi_{3} + \xi_{4}) - m(\xi_{2})| \lesssim \frac{|\xi_{3} + \xi_{4}|}{|\xi_{2}|}.
\end{equation}

\noindent Moreover, $(\ref{4.7})$ implies

\begin{equation}\label{4.9}
I((P_{> \frac{N}{8}} u)(P_{\leq \frac{N}{8}} u)^{2}) - (I P_{> \frac{N}{8}} u)(P_{\leq \frac{N}{8}} u)^{2}
\end{equation}

\noindent has a Fourier transform supported on $|\xi| \geq \frac{N}{8}$. Then by $(\ref{4.8})$, $E(Iu(t)) \leq 1$, and Theorem $\ref{t3.1}$,

\begin{equation}\label{4.10}
\int_{J} \langle i \nabla Iu, \nabla ((I P_{> \frac{N}{8}} u)(P_{\leq \frac{N}{8}} u)^{2} - I((P_{> \frac{N}{8}} u)(P_{\leq \frac{N}{8}} u)^{2})) \rangle dt
\end{equation}

\begin{equation}\label{4.11}
 \lesssim \frac{1}{N} \| \nabla I P_{> \frac{N}{8}} u \|_{L_{t}^{2} L_{x}^{6}(J \times \mathbf{R}^{3})}^{2} \| \nabla Iu \|_{L_{t}^{\infty} L_{x}^{2}(J \times \mathbf{R}^{3})} \| Iu \|_{L_{t}^{\infty} L_{x}^{6}(J \times \mathbf{R}^{3})} \lesssim \frac{1}{N}.
\end{equation}

\noindent Also since $E(Iu(t)) \leq 1$,

\begin{equation}\label{4.12}
 \int_{J} \langle i \nabla Iu, \nabla ((I P_{> \frac{N}{8}} u)^{2} (P_{\leq \frac{N}{8}} u) - I((P_{> \frac{N}{8}} u)^{2} (P_{\leq \frac{N}{8}} u))) \rangle dt
\end{equation}

\begin{equation}\label{4.13}
 \lesssim \| \nabla Iu \|_{L_{t}^{\infty} L_{x}^{2}} \| \nabla I P_{> \frac{N}{8}} u \|_{L_{t}^{2} L_{x}^{6}} \| I P_{> \frac{N}{8}} u \|_{L_{t}^{2} L_{x}^{6}} \| P_{\leq \frac{N}{8}} u \|_{L_{t}^{\infty} L_{x}^{6}} \lesssim \frac{1}{N}.
\end{equation}

\noindent Finally, by $(\ref{3.10})$ and $(\ref{3.11})$,

\begin{equation}\label{4.14}
\int_{J} \langle i \nabla Iu, \nabla ((I P_{> \frac{N}{8}} u)^{3} - I((P_{> \frac{N}{8}} u)^{3})) \rangle dt
\end{equation}

\begin{equation}\label{4.15}
 \lesssim \| \nabla Iu \|_{L_{t}^{\infty} L_{x}^{2}} \| \nabla I P_{> \frac{N}{8}} u \|_{L_{t}^{2} L_{x}^{6}} \| P_{> \frac{N}{8}} u \|_{L_{t}^{4} L_{x}^{6}}^{2} \lesssim \frac{1}{N}.
\end{equation}

\noindent This takes care of the first term in $(\ref{4.4})$. Now consider the term

\begin{equation}\label{4.16}
\int_{J} \langle I(|u|^{2} u), I(|u|^{2} u) - |Iu|^{2} (Iu) \rangle dt.
\end{equation}

\noindent $(\ref{4.6})$ and $(\ref{4.7})$ imply that this six-linear term must have at least two $P_{> \frac{N}{8}} u$ terms. By the Sobolev embedding theorem, Bernstein's inequality, and $(\ref{2.2})$,

\begin{equation}\label{4.17}
 \| I((P_{> \frac{N}{8}} u)^{3}) \|_{L_{t,x}^{2}} \lesssim \| \nabla I P_{> \frac{N}{8}} u \|_{L_{t}^{2} L_{x}^{6}} \| P_{> \frac{N}{8}} u \|_{L_{t}^{\infty} L_{x}^{3}}^{2} \lesssim \frac{1}{N}.
\end{equation}

\noindent Therefore,

\begin{equation}\label{4.18}
 \int_{J} \langle I((P_{> \frac{N}{8}} u)^{3}), I((P_{> \frac{N}{8}} u)^{3}) + (I P_{> \frac{N}{8}} u)^{3} \rangle dt \lesssim \frac{1}{N^{2}}.
\end{equation}

\noindent Next,

\begin{equation}\label{4.19}
\aligned
 \int_{J} \langle I((P_{> \frac{N}{8}} u)^{3}), (P_{> \frac{N}{8}} u)^{2} (P_{\leq \frac{N}{8}} u) \rangle dt \\ \lesssim \|  I (P_{> \frac{N}{8}} u)^{3} \|_{L_{t, x}^{2}} \| P_{> \frac{N}{8}} u \|_{L_{t}^{2} L_{x}^{6}} \| P_{\leq \frac{N}{8}} u \|_{L_{t,x}^{\infty}} \| P_{> \frac{N}{8}} u \|_{L_{t}^{\infty} L_{x}^{3}} \lesssim \frac{1}{N^{2}}.
\endaligned
\end{equation}

\noindent Finally,

\begin{equation}\label{4.20}
\int_{J} \int (P_{> \frac{N}{8}} u)^{2} (P_{\leq \frac{N}{8}} u)^{2} u^{2} dx dt
\end{equation}

\begin{equation}\label{4.21}
 \lesssim \| P_{> \frac{N}{8}} u \|_{L_{t}^{2} L_{x}^{6}}^{2} \| P_{\leq \frac{N}{8}} u \|_{L_{t}^{\infty} L_{x}^{6}}^{4} + \| P_{> \frac{N}{8}} u \|_{L_{t}^{4} L_{x}^{6}}^{4} \| P_{\leq \frac{N}{8}} u \|_{L_{t}^{\infty} L_{x}^{6}} \lesssim \frac{1}{N^{2}}.
\end{equation}

\noindent This proves Lemma $\ref{l4.1}$. $\Box$\vspace{5mm}

\noindent Rescaling back, we have proved

\begin{equation}\label{4.22}
\| u(t) \|_{L_{x}^{2}(\mathbf{R}^{3})} = \| u(0) \|_{L_{x}^{2}(\mathbf{R}^{3})},
\end{equation}

\noindent and

\begin{equation}\label{4.23}
\| u(t) \|_{\dot{H}^{s}(\mathbf{R}^{3})} \lesssim \| u(0) \|_{L^{2}(\mathbf{R}^{3})} + N^{\frac{1 - s}{2s - 1}} \| u(0) \|_{\dot{H}^{s}(\mathbf{R}^{3})}.
\end{equation}

\noindent Therefore, by $(\ref{3.24})$,

\begin{equation}\label{4.24}
\| u(t) \|_{H^{s}(\mathbf{R}^{3})} \lesssim C(s, \| u_{0} \|_{H^{s}(\mathbf{R}^{3})}) \| u_{0} \|_{H^{s}(\mathbf{R}^{3})},
\end{equation}

\noindent where $C$ behaves like $e^{C_{1} \frac{1 - s}{2s - 1}}$ for some constant $C_{1}$ as $s \searrow \frac{1}{2}$. This completes the proof of Theorem $\ref{t1.4}$ since $(\ref{4.24})$ gives a bound on $\| u \|_{L_{t,x}^{4}}$ by Theorem $\ref{t2.2}$. Interpolating this with the uniform bound on $\| u(t) \|_{H^{s}}$ implies a bound on $L_{t}^{p} L_{x}^{q}$, where $(p, q)$ is a $\frac{1}{2}$-admissible pair, that is $\frac{2}{p} = 3(\frac{1}{2} - \frac{1}{q} - \frac{1}{6})$. Since $s > \frac{1}{2}$, $p < \infty$. Partitioning $\mathbf{R}$ into finitely many pieces with $\| u \|_{L_{t}^{p} L_{x}^{q}(J_{l} \times \mathbf{R}^{3})} < \epsilon$ and making a perturbation argument, $\| u \|_{L_{t,x}^{5}(\mathbf{R} \times \mathbf{R}^{3})} < \infty$, which implies scattering. $\Box$

\section{A computed mass-critical bound}
In two dimensions the cubic problem

\begin{equation}\label{10.0}
i u_{t} + \Delta u = \mu |u|^{2} u, \hspace{5mm} u(0, x) = u_{0}, \hspace{5mm} \mu = \pm 1,
\end{equation}

\noindent is mass-critical (see $(\ref{1.4})$). $\mu = +1$ is the defocusing case and $\mu = -1$ is the focusing case.\vspace{5mm}

\noindent \cite{KTV} proved that $(\ref{10.0})$ was globally well-posed and scattering in the defocusing case ($\mu = 1$) and in the focusing case ($\mu = -1$) with mass less than the mass of the ground state. This result was extended to the nonradial case by \cite{D4}. However, \cite{KTV} and \cite{D4} did not compute an explicit bound, which we will do here for initial data lying in $\dot{H}^{s} \cap |x|^{s} L^{2} \subset L^{2}$.

\begin{theorem}\label{t10.1}
If $u_{0}$ is a radially symmetric function with $u_{0} \in \dot{H}^{s}(\mathbf{R}^{2})$, $s > 0$, then the defocusing initial value problem

\begin{equation}\label{10.1}
i u_{t} + \Delta u = |u|^{2} u, \hspace{5mm} u(0, x) = u_{0},
\end{equation}

\noindent has a solution on $[0, 1]$ with

\begin{equation}\label{10.2}
\| u \|_{L_{t,x}^{4}([0, 1] \times \mathbf{R}^{2})}^{4} \lesssim_{s} \| u_{0} \|_{\dot{H}^{s}(\mathbf{R}^{2})}^{\frac{8(1 - s)}{s} + 2} (1 + \| u_{0} \|_{L^{2}})^{\frac{4(1 - s)}{s} + 2}.
\end{equation}

\noindent The focusing initial value problem

\begin{equation}\label{10.1.1}
i u_{t} + \Delta u = -|u|^{2} u, \hspace{5mm} u(0, x) = u_{0},
\end{equation}

\noindent has a solution on $[0, 1]$ with

\begin{equation}\label{10.2.1}
\| u \|_{L_{t,x}^{4}([0, 1] \times \mathbf{R}^{2})}^{4} \lesssim_{s} \| u_{0} \|_{\dot{H}^{s}(\mathbf{R}^{2})}^{\frac{8(1 - s)}{s} + 2} (1 + \| u_{0} \|_{L^{2}})^{\frac{4(1 - s)}{s} + 2} (1 - \frac{\| u_{0} \|_{L^{2}}^{2}}{\| Q \|_{L^{2}}^{2}})^{-\frac{1}{s}},
\end{equation}

\noindent where $Q$ is the soliton for $(\ref{10.1.1})$, that is $Q$ solves the elliptic problem

\begin{equation}
\Delta Q + Q^{3} = Q.
\end{equation}

\end{theorem}

\noindent This gives a scattering result.

\begin{corollary}\label{c10.2}
The initial value problem $(\ref{10.1})$ is globally well-posed and scattering for initial data lying in $\dot{H}^{s}(\mathbf{R}^{2}) \cap |x|^{s} L^{2}(\mathbf{R}^{2})$.
\end{corollary}

\noindent \emph{Proof of corollary:} By time reversal symmetry it suffices to prove

\begin{equation}\label{10.3}
\| u \|_{L_{t,x}^{4}([0, \infty) \times \mathbf{R}^{2})} < \infty.
\end{equation}

\noindent Rescale so that $\| u_{0} \|_{\dot{H}^{s}} = \| |x|^{s} u_{0} \|_{L^{2}}$. Shift $t = 0$ to $t = 1$ and then make the pseudoconformal transformation, for $t > 0$,

\begin{equation}\label{10.4}
v(t, x) = \frac{1}{t} e^{i \frac{|x|^{2}}{4t}} u(\frac{-1}{t}, \frac{x}{t}).
\end{equation}

\noindent Then $v$ also solves $(\ref{10.1})$ with initial data

\begin{equation}\label{10.4.1}
\| v(-1, x) \|_{\dot{H}^{s}(\mathbf{R}^{2})} \lesssim \| u_{0} \|_{\dot{H}^{s}(\mathbf{R}^{2})} + \| |x|^{s} u_{0} \|_{L^{2}(\mathbf{R}^{2})}.
\end{equation}

\noindent Then by Theorem $\ref{t10.1}$,

\begin{equation}\label{10.5}
\| v \|_{L_{t,x}^{4}([-1, 0] \times \mathbf{R}^{2})} < \infty.
\end{equation}

\noindent It is easy to verify by direct computation that

\begin{equation}\label{10.5.1}
\| v \|_{L_{t,x}^{4}([-1, 0] \times \mathbf{R}^{2})} = \| u \|_{L_{t,x}^{4}([1, \infty) \times \mathbf{R}^{2})}.
\end{equation}

\noindent Then shifting $t = 1$ back to $t = 0$ gives $(\ref{10.3})$. $\Box$\vspace{5mm}

\noindent \emph{Proof of theorem $\ref{t10.1}$:} Without loss of generality suppose that $\| u_{0} \|_{\dot{H}^{s}} \gtrsim 1$. Otherwise, $(\ref{10.2})$ could be proved by a small data argument. Next choose

\begin{equation}\label{10.5.2}
\lambda \sim_{s} \| u_{0} \|_{\dot{H}^{s}}^{-\frac{4(1 - s)}{s} - 1} (1 + \| u_{0} \|_{L^{2}})^{-\frac{2(1 - s)}{s} - 1}.
\end{equation}

\noindent Then after rescaling by $(\ref{1.3})$,

\begin{equation}\label{10.5.3}
E(Iu(0)) \lesssim N^{2(1 - s)} \| u_{0} \|_{\dot{H}^{s}}^{-8(1 - s)} (1 + \| u_{0} \|_{L^{2}})^{-4(1 - s)}.
\end{equation}

\noindent Then if we choose 

\begin{equation}\label{10.5.3.1}
N \sim_{s} \| u_{0} \|_{\dot{H}^{s}}^{4} (1 + \| u_{0} \|_{L^{2}})^{2},
\end{equation}

\begin{equation}\label{10.5.4}
E(Iu(0)) \leq \frac{1}{2}.
\end{equation}

\noindent Then to prove Theorem $\ref{t10.1}$ it suffices to prove

\begin{equation}\label{10.6}
E(Iu(t)) \leq 1, \hspace{5mm} \text{for all} \hspace{5mm} t \in [0, \lambda^{-2}].
\end{equation}

\noindent As in the cubic problem in three dimensions, this result will be proved using a long-time Strichartz estimate. First observe that if $J$ is an interval, $J \subset [0, \lambda^{-2}]$, and $E(Iu(t)) \leq 1$ for all $t \in J$, then

\begin{equation}\label{10.7}
\| Iu \|_{L_{t,x}^{4}(J \times \mathbf{R}^{2})}^{4} \lesssim \lambda^{-2}.
\end{equation}

\noindent Bernstein's inequality and $E(Iu(t)) \leq 1$ on $J$ implies that $\| (1 - I) u \|_{L_{t}^{\infty} L_{x}^{2}(J \times \mathbf{R}^{2})} \lesssim N^{-1}$, so then by standard perturbative arguments, if $\| Iu \|_{L_{t,x}^{4}(I_{j} \times \mathbf{R}^{2})} \leq \epsilon$, then for $N$ sufficiently large,

\begin{equation}\label{10.8}
\| u \|_{L_{t,x}^{4}(I_{j} \times \mathbf{R}^{2})} \leq 2 \epsilon.
\end{equation}

\noindent Therefore,

\begin{equation}\label{10.9}
\| u \|_{L_{t,x}^{4}(J \times \mathbf{R}^{2})}^{4} \lesssim \lambda^{-2},
\end{equation}

\noindent and thus since $E(Iu(t)) \leq 1$ on $J$,

\begin{equation}\label{10.10}
\| \nabla I u \|_{U_{\Delta}^{2}(J \times \mathbf{R}^{2})} \lesssim \lambda^{-1}.
\end{equation}

\begin{theorem}[Long time Strichartz estimate]\label{t10.3}
If $E(Iu(t)) \leq 1$ on $J$, then

\begin{equation}\label{10.10.1}
\| \nabla I P_{> \frac{N}{8}} u \|_{U_{\Delta}^{2}(J \times \mathbf{R}^{2})} \lesssim 1.
\end{equation}
\end{theorem}

\noindent \emph{Proof:} Again by $(\ref{7.25})$, for any $M \leq N$,

\begin{equation}\label{10.10.2}
\aligned
\| \nabla I P_{> M} u \|_{U_{\Delta}^{2}(J \times \mathbf{R}^{2})} \lesssim \| \nabla I P_{> \frac{M}{8}} u(0) \|_{L^{2}} + \| \nabla I ((P_{> \frac{M}{8}} u)^{2} u) \|_{L_{t}^{2-} L_{x}^{1+}(J \times \mathbf{R}^{2})} \\ + \| (\nabla I P_{\leq \frac{M}{8}} u) (P_{> \frac{M}{8}} u)(P_{\leq \frac{M}{8}} u) \|_{L_{t}^{2-} L_{x}^{1+}(J \times \mathbf{R}^{2})} + \frac{1}{M^{1 - \frac{1}{q}}} \| (\nabla I P_{> \frac{M}{8}} u)(P_{\leq \frac{M}{8}} u)^{2} \|_{X_{R}}.
\endaligned
\end{equation}
\textbf{Remark:} Here we will use the $+$ and $-$ notation instead of $L_{t}^{2- \epsilon} L_{x}^{1 + \delta(\epsilon)}$, and will not explicitly compute the $\epsilon$ dependence of the exponents. The interested reader could use the analysis in section four as a template, since the computations are quite similar. The important fact is that the powers will be bounded by a constant times $\epsilon > 0$.\medskip

First, since $E(Iu(t)) \leq 1$ for all $t \in J$, $\| \nabla I P_{> \frac{M}{8}} u(0) \|_{L^{2}} \lesssim 1$. Next, by $(\ref{10.9})$ and Bernstein's inequality,

\begin{equation}\label{10.11}
\| \nabla I ((P_{> \frac{M}{8}} u)^{2} u) \|_{L_{t}^{2-} L_{x}^{1+}(J \times \mathbf{R}^{2})} \lesssim \frac{1}{M \lambda^{0+}} \| \nabla I P_{> \frac{M}{8}} u \|_{U_{\Delta}^{2}(J \times \mathbf{R}^{2})},
\end{equation}

\noindent and

\begin{equation}\label{10.12}
\| (\nabla I P_{\leq \frac{M}{8}} u)(P_{> \frac{M}{8}} u)(P_{\leq \frac{M}{8}} u) \|_{L_{t}^{2-} L_{x}^{1+}(J \times \mathbf{R}^{2})} \lesssim \frac{1}{M \lambda^{0+}} \| \nabla I P_{> \frac{M}{8}} u \|_{U_{\Delta}^{2}(J \times \mathbf{R}^{2})}.
\end{equation}

\noindent Finally, by the fundamental theorem of calculus and $E(Iu(t)) \leq 1$,

\begin{equation}\label{10.13}
|x| |P_{\leq \frac{M}{8}} u|^{2} \leq \int_{|x|}^{\infty} r \partial_{r} (|P_{\leq \frac{M}{8}} u|^{2}) dr \lesssim \| \nabla I u \|_{L^{2}} \| Iu \|_{L^{2}} \lesssim \| u_{0} \|_{L^{2}}.
\end{equation}

\noindent Then for any $j \geq 0$ and $R$, by $(\ref{10.10})$,

\begin{equation}\label{10.14}
\aligned
\frac{R^{1 - \frac{1}{q}} 2^{j(1 - \frac{1}{q})}}{M^{1 - \frac{1}{q}}} \| \chi(\frac{x}{2^{j} R}) (\nabla I P_{> \frac{M}{8}} u)(P_{\leq \frac{M}{8}} u)^{2} \|_{L_{t}^{q} L_{x}^{2}(J \times \mathbf{R}^{2})} \\ \lesssim \frac{\lambda^{\frac{1}{2} - \frac{1}{q}}}{M^{\frac{3}{2} - \frac{1}{q}}} \| u_{0} \|_{L^{2}} \| \nabla I P_{> \frac{M}{8}} u \|_{U_{\Delta}^{2}(J \times \mathbf{R}^{2})}.
\endaligned
\end{equation}

\noindent Also by $(\ref{10.13})$ and $J \subset [0, \lambda^{-2}]$,

\begin{equation}\label{10.15}
\aligned
R^{1 - \frac{1}{q}} 2^{j(1 - \frac{1}{q})} \| \chi(\frac{x}{2^{j} R}) (\nabla I P_{> \frac{M}{8}} u)(P_{\leq \frac{M}{8}} u)^{2} \|_{L_{t}^{q} L_{x}^{2}(J \times \mathbf{R}^{2})} \\ \lesssim R^{\frac{1}{2} - \frac{1}{q}} 2^{j(\frac{1}{2} - \frac{1}{q})} \frac{\lambda^{1 - \frac{2}{q}}}{M^{\frac{1}{2}}} \| u_{0} \|_{L^{2}} \| \nabla I P_{> \frac{M}{8}} u \|_{U_{\Delta}^{2}(J \times \mathbf{R}^{2})}.
\endaligned
\end{equation}

\noindent Also by $(\ref{10.9})$ and $(\ref{10.10})$,

\begin{equation}\label{10.16}
R^{1 - \frac{1}{q}} \| \psi(\frac{x}{R}) (\nabla I P_{> \frac{M}{8}} u)(P_{\leq \frac{M}{8}} u)^{2} \|_{L_{t}^{q} L_{x}^{2}(J \times \mathbf{R}^{2})} \lesssim \frac{\lambda^{1 - \frac{2}{q}} R^{1 - \frac{1}{q}}}{M^{\frac{1}{2}}} \|  \nabla I P_{> \frac{M}{8}} u \|_{U_{\Delta}^{2}(J \times \mathbf{R}^{2})}.
\end{equation}

\noindent Then taking $R = 1$ and using $(\ref{10.15})$ and $(\ref{10.16})$ to sum over $j$, combined with the fact that $\lambda \sim N^{\frac{-(1 - s)}{s}}$, and taking $q$ arbitrarily close to $2$,

\begin{equation}\label{10.17}
\| \nabla I P_{> M} u \|_{U_{\Delta}^{2}(J \times \mathbf{R}^{2})} \lesssim 1 + \frac{N^{0+}}{M^{1-}} \| u_{0} \|_{L^{2}} \| \nabla I P_{> \frac{M}{8}} u \|_{U_{\Delta}^{2}(J \times \mathbf{R}^{2})}.
\end{equation}

\noindent Then making an induction on frequency argument, starting with $M = N^{3/4}$, implies that

\begin{equation}\label{10.18}
\| \nabla I P_{> \frac{N}{8}} u \|_{U_{\Delta}^{2}(J \times \mathbf{R}^{2})} \lesssim 1 + N^{-c \ln(N)} \| u_{0} \|_{\dot{H}^{s}}^{\frac{4(1 - s)}{s} + 1} (1 + \| u_{0} \|_{L^{2}})^{\frac{2(1 - s)}{s} + 1}.
\end{equation}

\noindent Then if $N$ is given by $(\ref{10.5.3.1})$, the proof is complete. $\Box$\vspace{5mm}

This gives a bound on the growth of $E(Iu(t))$.

\begin{theorem}\label{t10.3}
$E(Iu(t)) \leq 1$ for all $t \in [0, \lambda^{-2}]$.
\end{theorem}

\noindent \emph{Proof:} Since $E(Iu(0)) \leq \frac{1}{2}$, it remains to bound the time integral of $\frac{d}{dt} E(Iu(t))$. Much of the analysis in section five may be copied directly to this situation as well. However, there are some differences due to the difference in dimension, and thus there are different exponents due to different Sobolev embeddings. For example, instead of $(\ref{4.11})$, estimate

\begin{equation}\label{10.18}
\aligned
\| \nabla Iu \|_{L_{t}^{\infty} L_{x}^{2}(J \times \mathbf{R}^{2})} \| \nabla I P_{> \frac{N}{8}} u \|_{L_{t}^{2+} L_{x}^{\infty-}(J \times \mathbf{R}^{2})} \| P_{> \frac{N}{8}} u \|_{L_{t}^{2+} L_{x}^{\infty-}(J \times \mathbf{R}^{2})} \| P_{\leq \frac{N}{8}} u \|_{L_{t}^{\infty-} L_{x}^{2+}(J \times \mathbf{R}^{2})} \\ \lesssim \frac{1}{N} \frac{1}{\lambda^{0+}} \| u_{0} \|_{L^{2}}.
\endaligned
\end{equation}

\noindent This finishes the proof of Theorem $\ref{t10.1}$ in the defocusing case.\vspace{5mm}

For the focusing problem use the Gagliardo-Nirenberg inequality (see \cite{Nirenberg}),

\begin{equation}\label{10.19}
\| u \|_{L_{x}^{4}(\mathbf{R}^{2})}^{4} \leq \frac{1}{2} \frac{\| u \|_{L_{x}^{2}(\mathbf{R}^{2})}^{2}}{\| Q \|_{L_{x}^{2}(\mathbf{R}^{2})}^{2}} \| \nabla u \|_{L_{x}^{2}(\mathbf{R}^{2})}^{2}.
\end{equation}

\noindent Therefore,

\begin{equation}\label{10.20}
\| \nabla I u \|_{L_{x}^{2}(\mathbf{R}^{2})}^{2} (1 - \frac{\| u_{0} \|_{L^{2}}^{2}}{\| Q \|_{L^{2}}^{2}}) \lesssim E(Iu),
\end{equation}

\noindent where in this case

\begin{equation}\label{10.21}
E(Iu(t)) = \frac{1}{2} \int |\nabla Iu(t,x)|^{2} dx - \frac{1}{4} \int |Iu(t,x)|^{4} dx.
\end{equation}

\noindent Then replace $(\ref{10.5.2})$ and $(\ref{10.5.3.1})$ with

\begin{equation}
\lambda \sim_{s} \| u_{0} \|_{\dot{H}^{s}}^{-\frac{4(1 - s)}{s} - 1} (1 + \| u_{0} \|_{L^{2}})^{-\frac{2(1 - s)}{s} - 1} (1 - \frac{\| u_{0} \|_{L^{2}}^{2}}{\| Q \|_{L^{2}}^{2}})^{\frac{1}{s}}
\end{equation}

\noindent and

\begin{equation}
N \sim_{s} \| u_{0} \|_{\dot{H}^{s}}^{4} (1 + \| u_{0} \|_{L^{2}})^{2} (1 - \frac{\| u_{0} \|_{L^{2}}^{2}}{\| Q \|_{L^{2}}^{2}})^{-1},
\end{equation}

\noindent respectively. Then proceed as in the defocusing case. $\Box$

\section{Induction on frequency in two dimensions}
Now turn to the two dimensional problem $(\ref{1.2})$ with $k > 1$, $k \in \mathbf{Z}$. Here the critical space is $\dot{H}^{s_{c}}$, $s_{c} = \frac{k - 1}{k}$. Once again take the $I$ operator as defined in $(\ref{2.2})$. Then,

\begin{equation}\label{5.2}
E(Iu(0)) \lesssim_{k, \| u_{0} \|_{H^{s}(\mathbf{R}^{2})}} N^{2(1 - s)}.
\end{equation}

\noindent Rescale with $\lambda \sim_{\| u_{0} \|_{H^{s}}, k} N^{\frac{1 - s}{s - s_{c}}}$ so that $E(Iu(0)) = \frac{1}{2}$. After rescaling $M(Iu(0)) \lesssim N^{\frac{1 - s}{2(s - s_{c})}}$. Suppose $J$ is an interval with $E(Iu(t)) \leq 1$ for all $t \in J$. Recalling $(\ref{2.20})$,

\begin{equation}\label{5.3}
\| u \|_{L_{t}^{4} L_{x}^{8}(J \times \mathbf{R}^{2})}^{4} \lesssim \| |\nabla|^{1/2} |u|^{2} \|_{L_{t,x}^{2}}^{2} \lesssim_{\| u(0) \|_{H^{s}(\mathbf{R}^{2})}, k} N^{s_{c} \cdot \frac{3(1 - s)}{s - s_{c}}}.
\end{equation}

\noindent Then, by Lemma $\ref{l2.3}$,

\begin{equation}\label{5.4}
\| \nabla Iu \|_{U_{\Delta}^{2}(J \times \mathbf{R}^{2})} \lesssim_{\| u(0) \|_{H^{s}(\mathbf{R}^{2})}, k} N^{s_{c} \frac{3(1 - s)}{2(s - s_{c})}}.
\end{equation}

Once make an induction on frequency argument to prove long time Strichartz estimates.

\begin{theorem}\label{t5.1}
Let $0 \in J$ be an interval such that $E(Iu(t)) \leq 1$ on $J$. Then for any $s > s_{c}$, there exists $N(s, k, \| u_{0} \|_{H^{s}}) < \infty$ such that

\begin{equation}\label{5.4.1}
\| \nabla I P_{> \frac{N}{8k}} u \|_{U_{\Delta}^{2}(J \times \mathbf{R}^{2})} \lesssim_{p} 1.
\end{equation}
\end{theorem}

\noindent \emph{Proof:} Again by $(\ref{7.25})$, if $M \leq N$,

\begin{equation}\label{5.5}
\aligned
\| \nabla I P_{> M} u(t) \|_{U_{\Delta}^{2}(J \times \mathbf{R}^{2})} \lesssim \| \nabla I P_{> M} u(0) \|_{L_{x}^{2}(\mathbf{R}^{2})} + \| \nabla I P_{> M} ((P_{> \frac{M}{8k}} u)^{2} u^{2k - 1}) \|_{L_{t}^{2-} L_{x}^{1+}(J \times \mathbf{R}^{2})} \\
+ \| (I P_{> \frac{M}{8k}} u) (\nabla P_{\leq \frac{M}{8k}} u)(P_{\leq \frac{M}{8k}} u)^{2k - 1} \|_{L_{t}^{2-} L_{x}^{1+}(J \times \mathbf{R}^{2})} + \frac{1}{M^{1 - \frac{1}{q}}} \| (\nabla I P_{> \frac{M}{8k}} u) (P_{\leq \frac{M}{8k}} u)^{2k} \|_{X_{R}}.
\endaligned
\end{equation}

\noindent Once again, since the nonlinearity is algebraic,

\begin{equation}\label{5.7}
P_{> M} (|u_{\leq \frac{M}{8k}}|^{2k} u_{\leq \frac{M}{8k}}) = 0.
\end{equation}

\noindent Once again it is also perfectly fine to not distinguish between $u$ and $\bar{u}$. Now again since the Fourier multiplier of $\nabla I$ is increasing as $|\xi| \nearrow \infty$,

\begin{equation}\label{5.8}
\| \nabla I((P_{> \frac{M}{8k}} u)^{2} u^{2k - 1}) \|_{L_{t}^{2-} L_{x}^{1+}(J \times \mathbf{R}^{2})}
\end{equation}

\begin{equation}\label{5.9}
\lesssim \| \nabla I P_{> \frac{M}{8k}} u \|_{L_{t}^{2} L_{x}^{\infty-}} \| P_{> \frac{M}{8k}} u \|_{L_{t}^{\infty-} L_{x}^{2+}} \| Iu \|_{L_{t}^{\infty-} L_{x}^{\infty}}^{k - 2} \| Iu \|_{L_{t}^{\infty} L_{x}^{2k + 2}}^{k + 1}
\end{equation}

\begin{equation}\label{5.10}
+ \| \nabla I P_{> \frac{M}{8k}} u \|_{L_{t}^{2+} L_{x}^{\infty-}} \| P_{> \frac{M}{8k}} u \|_{L_{t}^{\infty} L_{x}^{2k}} \| P_{> N} u \|_{L_{t}^{\infty-} L_{x}^{2k+}}^{2k - 1} 
\end{equation}

\begin{equation}\label{5.12}
+ \| P_{> \frac{M}{8k}} u \|_{L_{t, x}^{4}}^{2}  \| \nabla Iu \|_{L_{t}^{\infty-} L_{x}^{2+}} \| Iu \|_{L_{t}^{\infty-} L_{x}^{\infty}}^{2k - 1}
\end{equation}

\begin{equation}\label{5.11}
+ \| P_{> \frac{M}{8k}} u \|_{L_{t}^{4k} L_{x}^{4k}}^{2} \| \nabla Iu \|_{L_{t}^{\infty-} L_{x}^{2+}}  \| P_{> N} u \|_{L_{t}^{4k} L_{x}^{4k}}^{2k - 2}.
\end{equation}
\textbf{Remark:} Once again, we will use the $+$ and $-$ notation, rather than explicitly computing the $\epsilon$ dependence in the exponents.\medskip

\noindent Now by Lemma $\ref{l2.3}$ and $(\ref{5.4})$,

\begin{equation}\label{5.11.1}
\| Iu \|_{L_{t}^{\infty-} L_{x}^{\infty}} + \| \nabla Iu \|_{L_{t}^{\infty-} L_{x}^{2+}} \lesssim N^{+}.
\end{equation}

\noindent Also by interpolation and Bernstein's inequality,

\begin{equation}\label{5.11.2}
\| P_{> \frac{M}{8k}} u \|_{L_{t,x}^{4}}^{2} \lesssim \frac{1}{M} \| \nabla Iu \|_{L_{t}^{\infty} L_{x}^{2}} \| \nabla Iu \|_{U_{\Delta}^{2}},
\end{equation}

\noindent and

\begin{equation}\label{5.11.3}
\| P_{> \frac{M}{8k}} u \|_{L_{t,x}^{4k}}^{2k} \lesssim \frac{1}{M} \| \nabla Iu \|_{L_{t}^{\infty} L_{x}^{2}}^{2k - 1} \| \nabla Iu \|_{U_{\Delta}^{2}}.
\end{equation}

%\noindent Interpolating $(\ref{2.15})$ with $E(Iu(t)) \leq 1$ on $J$, combined with Bernstein's inequality %implies

%\begin{equation}\label{3.8}
%\| Iu \|_{L_{t}^{\infty-} L_{x}^{6+}(J \times \mathbf{R}^{3})} + \| \nabla Iu \|_{L_{t}^{\infty-} L_{x}^{2+}(J %\times \mathbf{R}^{3})}  \lesssim_{\| u_{0} \|_{H^{s}}, s} N^{+}.
%\end{equation}

%\noindent and

%\begin{equation}\label{3.9}
%\| P_{> N} u \|_{L_{t}^{\infty-} L_{x}^{3+}(J \times \mathbf{R}^{3})} \lesssim N^{-1/2+}.
%\end{equation}

%\noindent Interpolation, Sobolev embedding, Bernstein's inequality, $E(Iu(t)) \leq 1$, and $(\ref{2.2})$ imply %that

%\begin{equation}\label{3.10}
%\| P_{> \frac{M}{8}} u \|_{L_{t}^{4} L_{x}^{6}(J \times \mathbf{R}^{3})} \lesssim \| |\nabla|^{1/2} P_{> \frac{M}%{8}} u \|_{L_{t}^{2} L_{x}^{6}(J \times \mathbf{R}^{3})}^{1/2} \| |\nabla|^{1/2} P_{> \frac{M}{8}} u \|%_{L_{t}^{\infty} L_{x}^{2}(J \times \mathbf{R}^{3})}^{1/2}
%\end{equation}

%\begin{equation}\label{3.11}
%\lesssim M^{-1/2} \| \nabla I P_{> \frac{M}{8}} u \|_{L_{t}^{2} L_{x}^{6}(J \times \mathbf{R}^{3})}^{1/2}.
%\end{equation}

\noindent Making an argument almost identical to the estimates when $d = 3$,

\begin{equation}\label{5.13}
(\ref{5.8}) \lesssim_{k, \| u_{0} \|_{H^{s}(\mathbf{R}^{2})}} \frac{N^{+}}{M^{1-}} \| \nabla I P_{> \frac{M}{8k}} u \|_{U_{\Delta}^{2}(J \times \mathbf{R}^{2})}.
\end{equation}

\noindent Similarly, since $M \leq N$ and $E(Iu(t)) \leq 1$,

\begin{equation}\label{5.13.1}
\aligned
\| (I P_{> \frac{M}{8k}} u) (\nabla P_{\leq \frac{M}{8k}} u)(P_{\leq \frac{M}{8k}} u)^{2k - 1} \|_{L_{t}^{2-} L_{x}^{1+}(J \times \mathbf{R}^{2})} \\\lesssim \frac{1}{M} \| \nabla I P_{> \frac{M}{8k}} u \|_{L_{t}^{2+} L_{x}^{\infty-}} \| \nabla Iu \|_{L_{t}^{\infty-} L_{x}^{2+}} \| Iu \|_{L_{t}^{\infty-} L_{x}^{\infty}}^{k - 2} \| Iu \|_{L_{t}^{\infty} L_{x}^{2k + 2}}^{k + 1} \\ \lesssim \frac{N^{+}}{M} \| \nabla I P_{> \frac{M}{8k}} u \|_{U_{\Delta}^{2}(J \times \mathbf{R}^{2})}.
\endaligned
\end{equation}

\noindent Once again,

\begin{equation}\label{5.14}
\frac{1}{M^{1 - \frac{1}{q}}} \| (\nabla I P_{> \frac{M}{8k}} u)(P_{\leq \frac{M}{8k}} u)^{2k} \|_{X_{R}(J \times \mathbf{R}^{2})}
\end{equation}

\noindent is estimated by the local smoothing estimate

%\begin{equation}\label{5.15}
 %\| |\nabla|^{1/2} e^{it \Delta} u_{0} \|_{L_{t,x}^{2}(\mathbf{R} \times \{ |x| \leq R \})} \lesssim R^{1/2} \| u_{0} \|_{L^{2}},
%\end{equation}

%\noindent as well as,

\begin{equation}\label{5.16}
\| \nabla I P_{> \frac{M}{8k}} u \|_{L_{t,x}^{2}(J \times \{ x : |x| \leq R \})} \lesssim \frac{R^{1/2}}{M^{1/2}} \| \nabla I P_{> \frac{M}{8k}} u \|_{U_{\Delta}^{2}(J \times \mathbf{R}^{2})}.
\end{equation}

%\noindent Let $\chi \in C_{0}^{\infty}(\mathbf{R}^{2})$, $\chi \equiv 1$ on $|x| \leq 1$, $\chi$ supported on $|x| \leq 2$. Take $\| v \|_{V_{\Delta}^{2}(J \times \mathbf{R}^{2})} = 1$, $\hat{v}(\xi)$ supported on $|\xi| \geq M$. Then when $k \geq 2$,

%\begin{equation}\label{5.17}
%\aligned
%\int_{J} \langle v, \chi^{2} \nabla I((P_{> c(k)M} u)(P_{\leq c(k)M} u)^{2k}) \rangle dt \\ \lesssim \| \chi v \|_{L_{t}^{2+} L_{x}^{\infty-}(\mathbf{R} \times \mathbf{R}^{2})} \| \chi \nabla I P_{> c(k)M} \|_{L_{t,x}^{2}(\mathbf{R} \times \mathbf{R}^{2})} \| Iu \|_{L_{t}^{\infty-} L_{x}^{4k+}(J \times \mathbf{R}^{2})}^{2k} \\ \lesssim N^{+} M^{-1/2} \| \nabla I P_{> c(k)M} u \|_{U_{\Delta}^{2}(J \times \mathbf{R}^{2})}.
%\endaligned
%\end{equation}

%\noindent Now take the cutoff supported on the annulus $|x| \sim 2^{j}$, $\psi_{j}(x) = \chi^{2}(2^{-j} x) - \chi^{2}(2^{-j + 1} x)$, where $j \geq 0$. When $k = 2$,

%\noindent When $k = 2$,

%\begin{equation}\label{5.18}
%\| |\nabla|^{1/2} (Iu)^{2} \|_{L_{x}^{2}(\mathbf{R}^{2})} \lesssim \| |\nabla|^{1/2} Iu \|_{L_{x}^{3}(\mathbf{R}^{2})} \| Iu \|_{L_{x}^{6}(\mathbf{R}^{2})} \lesssim E(Iu)^{2/3} \leq 1.
%\end{equation}

\noindent Then by the fundamental theorem of calculus,

\begin{equation}\label{5.19}
|Iu(x)|^{2k} \lesssim \frac{1}{|x|} \int_{|x|}^{\infty} r \partial_{r}(|Iu(r)|^{2k}) dr \lesssim \frac{1}{|x|} \| \nabla Iu \|_{L^{2}} \| Iu \|_{L_{x}^{4k - 2}}^{2k - 1} \lesssim \frac{1}{|x|}.
\end{equation}

\noindent The last inequality follows from the fact that $E(Iu(t)) \leq 1$ along with the interpolation (for $k > 1$)

\begin{equation}\label{5.19.1}
\| Iu \|_{L_{x}^{4k - 2}} \lesssim \| \nabla Iu \|_{L^{2}}^{\theta} \| Iu \|_{L^{2k + 2}}^{1 - \theta}.
\end{equation}

\noindent It is not particularly important what $\theta$ is. Then by $(\ref{5.19})$ and $(\ref{5.11.1})$, for any $j \geq 0$,

\begin{equation}\label{5.20}
\frac{2^{j(1 - \frac{1}{q})} R^{1 - \frac{1}{q}}}{M^{1 - \frac{1}{q}}} \| \chi(\frac{2^{-j} x}{R}) (P_{< \frac{M}{8k}} u)^{2k} (\nabla I P_{> \frac{M}{8k}} u) \|_{L_{t}^{q} L_{x}^{2}} \lesssim \frac{N^{+}}{M^{\frac{3}{2} - \frac{1}{q}}} \| \nabla I P_{> \frac{M}{8k}} u \|_{U_{\Delta}^{2}(J \times \mathbf{R}^{2})}.
\end{equation}

%\noindent \textbf{Remark:} For functions in $\dot{H}^{1/2}$, the radial Sobolev embedding theorem actually only holds for each Littlewood - Paley projection. However, it is possible to utilize the fact that $Iu$ has finite mass, $\| Iu \|_{L^{2}} \lesssim N^{\frac{1 - s}{2(s - s_{c})}}$ at very low frequencies, and then use $(\ref{5.20})$ for frequencies above the low frequency threshold at the price of a factor on the order of $\ln(N)$.\vspace{5mm}

\noindent Then for $j$ very large, for $k \geq 2$,

\begin{equation}\label{5.21}
\aligned
\sum_{j \geq J} \frac{2^{j(1 - \frac{1}{q})} R^{1 - \frac{1}{q}}}{M^{1 - \frac{1}{q}}} \| \chi(\frac{2^{-j} x}{R}) (P_{< \frac{M}{8k}} u)^{2k}(\nabla I P_{> \frac{M}{8k}} u) \|_{L_{t}^{q} L_{x}^{2}(J \times \mathbf{R}^{2})} \\ \lesssim \sum_{j \geq J} \frac{2^{j(\frac{1}{q} - \frac{1}{2})} R^{\frac{1}{q} - \frac{1}{2}}}{M^{\frac{3}{2} - \frac{1}{q}}} \| Iu \|_{L_{t}^{\infty} L^{2}} \| Iu \|_{L_{t}^{\infty} \dot{H}^{1}} \| \nabla I P_{> \frac{M}{8k}} u \|_{U_{\Delta}^{2}(J \times \mathbf{R}^{2})} \| Iu \|_{L_{t}^{\infty-} L_{x}^{\infty}}^{2k - 2} \\
\lesssim \sum_{j \geq J} \frac{2^{j(\frac{1}{q} - \frac{1}{2})} R^{\frac{1}{q} - \frac{1}{2}}}{M^{\frac{3}{2} - \frac{1}{q}}} N^{\frac{1 - s}{2(s - s_{c})}} N^{+}  \| \nabla I P_{> \frac{M}{8k}} u \|_{U_{\Delta}^{2}(J \times \mathbf{R}^{2})}.
\endaligned
\end{equation}

\noindent Then taking $J(N, s, s_{c}, k, R)$ sufficiently large and $q$ arbitrarily close to $2$,

\begin{equation}\label{5.22}
\lesssim \frac{N^{+}}{M^{1-}} \| \nabla I P_{> \frac{M}{8k}} u \|_{U_{\Delta}^{2}(J \times \mathbf{R}^{2})}.
\end{equation}

\noindent Finally taking $R(N)$ sufficiently small, by the Sobolev embedding theorem, and $(\ref{5.11.1})$,

%\noindent Finally, $E(Iu(t)) \leq 1$ combined with $\| u(t) \|_{L_{x}^{2}(\mathbf{R}^{3})} \lesssim N^{\frac{1 %- s}{2s - 1}}$ implies that

%\begin{equation}\label{3.22}
%\aligned
%\int_{J} \langle v, (1 - \chi(2^{-j} x)) (\nabla I P_{> \frac{M}{8}} u) (P_{\leq \frac{M}{8}} u)^{2} \rangle dt \\ %\lesssim \| v \|_{L_{t}^{4} L_{x}^{3}(J \times \mathbf{R}^{3})} \| \nabla P_{> \frac{M}{8}} Iu \|_{L_{t}^{2} %L_{x}^{6}(J \times \mathbf{R}^{3})} \| Iu \|_{L_{t,x}^{4}(J \times \mathbf{R}^{3})} \\ \times \| Iu \|%_{L_{t}^{\infty} L_{x}^{2}(J \times \mathbf{R}^{3})} ^{1/2} \| (1 - \chi(2^{-j} x)) Iu \|_{L_{t,x}^{\infty}(J \times %\mathbf{R}^{3})}^{1/2} \\ \lesssim 2^{-j/2} N^{\frac{1 - s}{s - 1/2}} \| \nabla I P_{> \frac{M}{2}} \|%_{U_{\Delta}^{2}(J \times \mathbf{R}^{3})}.
%\endaligned
%\end{equation}

%\noindent Combining $(\ref{3.12})$, $(\ref{3.16})$, $(\ref{3.21})$, and $(\ref{3.22})$, for $j \sim \ln (N)$ %sufficiently large,

%\noindent Therefore,

\begin{equation}\label{5.23}
\frac{R^{1 - \frac{1}{q}}}{M^{1 - \frac{1}{q}}} \| \psi(\frac{x}{R}) (\nabla I P_{> \frac{M}{8k}} u)(Iu)^{2k} \|_{L_{t}^{q} L_{x}^{2}(J \times \mathbf{R}^{2})} \lesssim \frac{R^{1 - \frac{1}{q}}}{M^{\frac{3}{2} - \frac{1}{q}}} N^{+} \| \nabla I P_{> \frac{M}{8k}} u \|_{U_{\Delta}^{2}(J \times \mathbf{R}^{2})}.
\end{equation}

\noindent This time we starting the induction at $C(s, \| u_{0} \|_{H^{s}}, k) N^{3/4}$ for $C(s, \| u_{0} \|_{H^{s}}, k)$ sufficiently large,

\begin{equation}\label{5.24}
 \| \nabla I P_{> \frac{N}{8k}} u \|_{U_{\Delta}^{2}(J \times \mathbf{R}^{2})} \lesssim_{\| u_{0} \|_{H^{s}}, s, k} 1 + N^{\frac{3(1 - s)}{4(s - s_{c})} s_{c}} N^{-\frac{c}{6} \ln(N)},
\end{equation}

\noindent for some constant $c > 0$.\vspace{5mm}

\noindent \textbf{Remark:} We could replace $\frac{c}{6}$ with $\frac{c}{q}$ for any $q > 4$. Therefore, choosing $N$ sufficiently large,

\begin{equation}\label{5.25}
 \| \nabla I P_{> \frac{N}{8k}} u \|_{U_{\Delta}^{2}(J \times \mathbf{R}^{3})} \lesssim_{\| u_{0} \|_{H^{s}}, k} 1.
\end{equation}

\section{Energy Increment}
\noindent To complete the proof of Theorem $\ref{t1.5}$ it remains to prove the usual bound on the growth of $E(Iu(t))$.

\begin{lemma}\label{l6.1}
If $J$ is an interval with $E(Iu(t)) \leq 1$ on $J$,

\begin{equation}\label{6.1}
\int_{J} |\frac{d}{dt} E(Iu(t))| dt \lesssim_{k} \frac{1}{N^{1-}}.
\end{equation}
\end{lemma}

%\noindent Recall that $(\ref{1.1})$ implies

%\begin{equation}\label{4.2}
%i Iu_{t} + \Delta Iu = |Iu|^{2} (Iu) + I(|u|^{2} u) - |Iu|^{2} (Iu).
%\end{equation}

\noindent \emph{Proof:} We compute

\begin{equation}\label{6.2}
\aligned
\frac{d}{dt} E(Iu(t)) = \langle Iu_{t}, |Iu|^{2k}(Iu) - I(|u|^{2k} u) \rangle \\
= -\langle i \nabla Iu, \nabla (|Iu|^{2k} (Iu) - I(|u|^{2k} u)) \rangle \\ - \langle i I(|u|^{2k} u), (|Iu|^{2k}(Iu) - I(|u|^{2k} u)) \rangle.
\endaligned
\end{equation}

\noindent Once again,

\begin{equation}\label{6.5}
 (I P_{\leq \frac{N}{8k}} u)^{2k + 1} - I((P_{\leq \frac{N}{8k}} u)^{2k + 1}) = 0.
\end{equation}

\noindent Also,

\begin{equation}\label{6.6}
\aligned
(I P_{> \frac{N}{8k}} u)(I P_{\leq \frac{N}{8k}} u)^{2k} - I((P_{> \frac{N}{8k}} u)(P_{\leq \frac{N}{8k}} u)^{2k}) \\ = (I P_{> \frac{N}{2}} u)(P_{\leq \frac{N}{8k}} u)^{2k} - I((P_{> \frac{N}{2}} u)(P_{\leq \frac{N}{8k}} u)^{2k}).
\endaligned
\end{equation}

%\noindent Again apply the fundamental theorem of calculus,

%\begin{equation}\label{6.7}
% |m(\xi_{2} + \xi_{3} + \xi_{4}) - m(\xi_{2})| \lesssim \frac{|\xi_{3} + \xi_{4}|}{|\xi_{2}|}.
%\end{equation}

%\noindent Moreover,

%\begin{equation}\label{4.9}
%I((P_{> \frac{N}{8}} u)(P_{\leq \frac{N}{8}} u)^{2}) - (I P_{> \frac{N}{8}} u)(P_{\leq \frac{N}{8}} u)^{2}
%\end{equation}

%\noindent has Fourier transform supported on $|\xi| \geq \frac{N}{8}$. Therefore, by Bernstein's inequality

\noindent As before in $(\ref{4.8})$,

\begin{equation}\label{6.8}
-\int_{J} \langle i \nabla Iu, \nabla ((I P_{> \frac{N}{8k}} u)(P_{\leq \frac{N}{8k}} u)^{2k} - I((P_{> \frac{N}{8k}} u)(P_{\leq \frac{N}{8k}} u)^{2k})) \rangle dt
\end{equation}

\begin{equation}\label{6.9}
 \lesssim \frac{1}{N} \| \nabla I P_{> \frac{N}{8k}} u \|_{L_{t}^{2+} L_{x}^{\infty-}(J \times \mathbf{R}^{2})}^{2} \| \nabla Iu \|_{L_{t}^{\infty-} L_{x}^{2+}(J \times \mathbf{R}^{2})} \| Iu \|_{L_{t}^{\infty} L_{x}^{(4k - 2)}(J \times \mathbf{R}^{2})}^{2k - 1} \lesssim \frac{1}{N^{1-}}.
\end{equation}

\noindent This follows from $(\ref{5.19.1})$, $(\ref{5.25})$ to estimate $\| \nabla I P_{> \frac{N}{8k}} u \|_{L_{t}^{2+} L_{x}^{\infty-}(J \times \mathbf{R}^{2})}$ and $(\ref{5.11.1})$ to estimate $\| \nabla Iu \|_{L_{t}^{\infty-} L_{x}^{2+}(J \times \mathbf{R}^{2})}$.\vspace{5mm}

\noindent Next, since $E(Iu(t)) \leq 1$,

\begin{equation}\label{6.10}
\int_{J} \langle i \nabla Iu, \nabla ((I P_{> \frac{N}{8k}} u)^{2} (P_{\leq \frac{N}{8k}} u)^{2k - 1} - I((P_{> \frac{N}{8k}} u)^{2} (P_{\leq \frac{N}{8k}} u)^{2k - 1})) \rangle dt
\end{equation}

\begin{equation}\label{6.11}
 \lesssim \| \nabla Iu \|_{L_{t}^{\infty-} L_{x}^{2+}} \| \nabla I P_{> \frac{N}{8k}} u \|_{L_{t}^{2+} L_{x}^{\infty-}} \| I P_{> \frac{N}{8k}} u \|_{L_{t}^{2+} L_{x}^{\infty-}} \| P_{\leq \frac{N}{8k}} u \|_{L_{t}^{\infty} L_{x}^{4k - 2}}^{2k - 1} \lesssim \frac{1}{N^{1-}}.
\end{equation}

\noindent Finally, we skip ahead to

\begin{equation}\label{6.12}
\int_{J} \langle i \nabla Iu, \nabla ((I P_{> \frac{N}{8k}} u)^{2k + 1} - I((P_{> \frac{N}{8k}} u)^{2k + 1})) \rangle dt
\end{equation}

\begin{equation}\label{6.13}
 \lesssim \| \nabla Iu \|_{L_{t}^{\infty-} L_{x}^{2+}} \| \nabla I P_{> \frac{N}{8k}} u \|_{L_{t}^{2+} L_{x}^{\infty-}} \| P_{> \frac{N}{8k}} u \|_{L_{t, x}^{4k}}^{2k} \lesssim \frac{1}{N^{1-}}.
\end{equation}

\noindent \textbf{Remark:} The other terms can be handled in a similar manner.\vspace{5mm}

%\noindent \textbf{Remark:} $L_{t,x}^{4k}$ is $\dot{H}^{\frac{k - 1}{k}}$ - critical. Interpolation shows that

%\begin{equation}\label{6.14}
%\| P_{> c(k)N} u \|_{L_{t,x}^{4k}}^{2k} \lesssim E(Iu(t))^{1/2 \cdot (2k - 1)-} \| \nabla I P_{> c(k)N} u \|_{L_{t}^{2+} L_{x}^{\infty-}}^{1+}.
%\end{equation}

\noindent Then by $(\ref{5.25})$, $E(Iu(t)) \leq 1$, we are done with the first term in $(\ref{6.2})$. Now we consider the term

\begin{equation}\label{6.15}
\int_{J} \langle I(u^{2k + 1}), I(u^{2k + 1}) - (Iu)^{2k + 1} \rangle dt.
\end{equation}

\noindent Once again this term must have at least two $P_{> \frac{N}{8k}} u$ terms. By the Sobolev embedding theorem, $E(Iu(t)) \leq 1$, Bernstein's inequality, $(\ref{5.11.1})$, and $(\ref{5.25})$,

\begin{equation}\label{6.16}
 \| I((P_{> \frac{N}{8k}} u)^{2k + 1}) \|_{L_{t,x}^{2}} \lesssim \| \nabla I P_{> \frac{N}{8k}} u \|_{L_{t}^{2+} L_{x}^{\infty-}} \| P_{> \frac{N}{8k}} u \|_{L_{t}^{\infty-} L_{x}^{2k+}}^{2k} \lesssim \frac{1}{N^{1-}}.
\end{equation}

\noindent Therefore,

\begin{equation}\label{6.17}
 \int_{J} \langle I((P_{> \frac{N}{8k}} u)^{2k + 1}), I((P_{> \frac{N}{8k}} u)^{2k + 1}) - (I P_{> \frac{N}{8k}} u)^{2k + 1} \rangle dt \lesssim \frac{1}{N^{2-}}.
\end{equation}

\noindent Next,

\begin{equation}\label{6.18}
\aligned
 \int_{J} \langle I((P_{> \frac{N}{8k}} u)^{2k + 1}), (P_{> \frac{N}{8k}} u)^{2k} (P_{\leq \frac{N}{8k}} u) \rangle dt \\ \lesssim \|  I (P_{> \frac{N}{8k}} u)^{2k + 1} \|_{L_{t, x}^{2}} \| P_{> \frac{N}{8k}} u \|_{L_{t, x}^{4k}}^{2k} \| P_{\leq \frac{N}{8k}} u \|_{L_{t,x}^{\infty}} \lesssim \frac{1}{N^{2-}}.
\endaligned
\end{equation}

\noindent Finally,

\begin{equation}\label{6.19}
\aligned
\int_{J} \int (P_{> \frac{N}{8k}} u)^{2} (P_{\leq \frac{N}{8k}} u)^{2} u^{4k - 2} dx dt \\ \lesssim \int_{J} \int (P_{> \frac{N}{8k}} u)^{4k + 2} dx dt + \int_{J} \int (P_{> \frac{N}{8k}} u)^{2}(P_{\leq \frac{N}{8k}} u)^{4k} dx dt.
\endaligned
\end{equation}

\noindent Interpolating the $L_{x}^{2k + 2}$ and $\dot{H}^{1}$ norms, since $E(Iu(t)) \leq 1$,

\begin{equation}\label{6.20}
 \| Iu \|_{L_{t}^{\infty} L_{x}^{4k}} \lesssim 1.
\end{equation}

\noindent This proves Lemma $\ref{l6.1}$. $\Box$\vspace{5mm}

\noindent Rescaling back, we have proved

\begin{equation}\label{6.21}
\| u(t) \|_{L_{x}^{2}(\mathbf{R}^{2})} = \| u(0) \|_{L_{x}^{2}(\mathbf{R}^{2})},
\end{equation}

\noindent and

\begin{equation}\label{6.22}
\| u(t) \|_{\dot{H}^{s}(\mathbf{R}^{2})} \lesssim \| u(0) \|_{L^{2}(\mathbf{R}^{2})} + N^{s_{c} \cdot \frac{1 - s}{s - s_{c}}} \| u(0) \|_{\dot{H}^{s}(\mathbf{R}^{2})}.
\end{equation}

\noindent Therefore,

\begin{equation}\label{6.23}
\| u(t) \|_{H^{s}(\mathbf{R}^{2})} \lesssim C(\| u_{0} \|_{H^{s}(\mathbf{R}^{2})}, k) \| u_{0} \|_{H^{s}(\mathbf{R}^{2})},
\end{equation}

\noindent where $C$ behaves like $e^{C_{1} \frac{s_{c} (1 - s)}{s - s_{c}}}$ for some constant $C_{1}$.

\nocite*
\bibliographystyle{plain}

\end{document}